%Version: S.Finashin 03-02-97  05:37:55pm
%Format: plain
\input amstex
\documentstyle{amsppt}

\def\Int{\mathop{\roman{Int}}\nolimits}
\def\Cl{\mathop{\roman{Cl}}\nolimits}
\def\tr{\mathop{\roman{tr}}\nolimits}

\def\codim{\mathop{\roman{codim}}\nolimits}
\def\conj{\mathop{\roman{conj}}\nolimits}
\def\rank{\mathop{\roman{rk}}\nolimits}
\def\sign{\mathop{\roman{sign}}\nolimits}
\def\lk{\mathop\roman{lk}\nolimits}
\def\d{\partial}
\def\ind{\mathop\roman{ind}\nolimits}
\def\Spin{\mathop\roman{Spin}\nolimits}

\def\Sing{\mathop\roman{Sing}\nolimits}
\def\A{\frak A}
%
%             Greek and Gothic letters
%
\def\a{\alpha}
\def\b{\beta}
\def\g{\gamma}
\def\G{\Gamma}
\def\D{\Delta}
\def\e{\varepsilon}
\def\s{\sigma}
\def\om{\omega}
%
%      Bbb
%
\def\C{{\Bbb C}}
\def\R{{\Bbb R}}
\def\Z{{\Bbb Z}}
\def\Q{{\Bbb Q}}

%
%
%       To get rid of AmSTeX's logo:
\def\nologo{\let\logo@\relax}
%
%                Special signs
%
\let\ge\geqslant
\let\le\leqslant
\def\la{\langle}
\def\ra{\rangle}
\def\til#1{\widetilde{#1}}
\def\Im{\roman{Im}}

\def\Cp#1{\C\roman P^{#1}}
\def\CP#1{\C\roman P^{#1}}
\def\Rp#1{\R\roman P^{#1}}
%
%                Abbreviations
%
\let\tm\proclaim
\let\endtm\endproclaim
\let\rk=\remark
\let\endrk=\endremark
\def\oo{\varnothing}

%%%%%%%%%%%%

\def\AG{A$\G$}
\def\ch{\roman{ch}}
\def\Lk{\mathop\roman{Lk}\nolimits}
\def\Lef{\mathop\roman{L}\nolimits}

\def\q{\roman q}

\def\re{\R}
\def\kapnew{\kappa}
\def\ee{\varkappa}

\def\ord{\mathop\roman{ord}\nolimits}

\def\Discr{\mathop\roman{Discr}\nolimits}
\def\bra{\mathop\roman{br}_{\roman{im}}\nolimits}
\def\top{\mathop\roman{top}\nolimits}
\def\Con{\mathop\roman{Con}\nolimits}
\def\res{\mathop\roman{res}\nolimits}
\def\QR{$\Q\roman{I}$}
\def\QQ{$\overline{\Q\roman{I}}\-$}
\def\sQ{$\Q\roman{I}^S$}
\def\QS{$\Q\roman{I}^S$}
\def\sQR{$\Q\roman{I}^S$}
\def\sQQ{$\overline{\Q\roman{I}}^S$}
\def\l{\lambda}
\def\B{\Theta}
\def\|{\,|\,}

\magnification\magstep1
\NoBlackBoxes
\NoRunningHeads

%\comment
\topmatter
\title
Complex intersection of real cycles in real algebraic varieties
%with inessential singularities
and generalized Arnold--Viro inequalities
\endtitle
\author
Sergey Finashin
\endauthor
\address
Middle East Technical University, Ankara, Turkey
\endaddress{}
\email
serge\,\@\,rorqual.cc.metu.edu.tr,
%finash\,\@\,pdmi.ras.ru
\endemail
%\date{}
\thanks
Partially supported by the Russian Foundation for Basic Researches,
grant 98-01-00832
\endthanks
%\keywords{}
%\abstract{}
\endtopmatter

\document

\heading
\S1. Introduction
\endheading
\subheading{1.1. The subject of the paper}
Considering a real algebraic variety,
we denote by $\C X$ its complex point set, by $\R X$ its
real part being the fixed point set of the complex conjugation,
$\conj\:\C X\to \C X$, and by $\overline X$ the quotient $\C X/\conj$.

Assuming that $\dim_{\C}\C X=d=2n$ and that $\C X$ is a rational
homology manifold, we obtain a $\Q$-valued quadratic form, $\psi$,
induced on $H_d(\R X)$ by the inclusion homomorphism from
the intersection form in $\C X$.
More generally, we can require that $\C X$ is a $\Q$-homology
manifold only in a neighborhood of $\R X$, or even more generally,
that $\overline X$ is a $\Q$-homology manifold in a
neighborhood of $\R X$ (in the latter case, we consider the quadratic form
on $H_d(\R X)$ induced from $\overline X$).
The form $\psi$ will be called {\it the complex intersection form} in $\R X$.
Its analysis gives some information about the topology of $\R X$.
For instance, in the case of a non-singular surface we obtain
the Arnold inequalities, which I tried to extend to the case of
singular varieties.

\subheading{1.2. The Petrovskii and Arnold inequalities}
Given a non-singular real algebraic curve $\C A\subset\Cp2$ of degree
$2k$, recall the Petrovskii inequalities
$$
\frac32k(k-1)\le p-n\le\frac32k(k-1)+1
\tag1-1
$$
where $p$ and $n$ denote the number of the components (called {\it ovals})
of $\R A$
lying inside even and, respectively, odd number of the other components.
This result is related to the 16th Hilbert problem and gives, for instance,
a negative answer to a particular question mentioned by  Hilbert: are there
non-singular real sextics, with the real part, $\R A$,
having 11 ovals
(11 is an upper bound for the number of ovals of a sextic
 provided by the Harnack theorem)
which bound disjoint topological discs in $\Rp2$.

To understand the nature of these inequalities, consider the double plane
$\pi\:X\to\Cp2$, branched along $\C A$. We can identify $X$ with one of
the real algebraic surfaces, $\C X^\pm$, defined by the equation
$f(x,y,z)=\pm t^2$ in a quasi-homogeneous complex projective $3$-space,
where $f$ is a real homogeneous polynomial of degree $2k$ defining $\C A$.
Note that $\R X^\pm$ is mapped by $\pi$ into the region
$\Rp2_\pm=\{[x:y:z]\in\Rp2\,|\,\pm f(x,y,z)\ge0\}$,
one of which is orientable, and the other is not.
Being free to vary the sign of $f$, we can assume that
the orientable region is $\Rp2_+$, then
$\chi(\R X^+)=2\chi(\Rp2_+)=2(p-n)$, $\chi(\R X^-)=2\chi(\Rp2_-)=2(n-p+1)$,
and the Petrovskii inequalities can be formulated as the estimates
$$
-\frac12\chi(\R X^\pm)\le\frac12h^{1,1}(X)-1
\tag1-2
$$
Note that the estimate (1-2), which belongs to Comessatti,
historically precedes the Petrovskii inequalities, although its
application to plane real curves in connection with
Hilbert's problem was found later.

To prove (1-2) it is enough to use the Riemann--Hurwitz formula in
combination with the formula for the signature of involution applied to
the branched coverings $\C X^\e\to \overline X^\e$ and $\C X^\e\to\Cp2$,
$\e=+$ or $\e=-$, which gives
$$
\aligned
2\chi(\overline X^\e)-\chi(\R X^\e)&=\chi(X)=2\chi(\Cp2)-
\chi(\C A)=4k^2-6k+6\\
2\s(\overline X^\e)-\la\R X^\e,\R X^\e\ra_X
&=\s(X)=2\s(\Cp2)-\la\C A,\C A\ra_X=2-2k^2,
\endaligned
\tag1-3
$$
where $\la\centerdot,\centerdot\ra_X$ stands for
the intersection indices in $X$.
Using that $b_1(X)=b_1(\overline X^\e)=0$
and $\la\R X^\e,\R X^\e\ra_X=-\chi(\R X^\e)$, since the tangent
and the normal bundle to $\R X^\e$ in $\C X^\e$ are anti-isomorphic,
we have
$$
\align
b_2^+(\overline X^\e)&=\frac12(b_2^+(X)-1)=p_g(X)=
\frac12(k-1)(k-2)\tag1-4\\
b_2^-(\overline X^\e)-\frac12\chi(\R X^\e)&=\frac12(b_2^-(X)-1)=
\frac12h^{1,1}(X)-1=\frac32k(k-1),
\tag1-5\endalign
$$
where (1-5) shows that (1-2) and the Petrovskii inequality
may be interpreted as $0\le b_2^-(\overline X^\e)$.
The estimates (1-4) and (1-5) (and their corollaries (1-1), (1-2))
can be further enforced if we take into account existence of certain
positive (or negative) square elements in $H_2(\overline X^\e)$.
There are known two such enforcements. One is due to
Comessatti, who used algebraic cycles to produce such elements,
the other is due to Arnold \cite{Ar}, who used the components of
$\R X^\e$ instead.
More precisely, Arnold's inequalities appear from the estimates
$$
\aligned
c^\pm&\le b_2^\pm(\overline X^\e)\\
c^\pm+c^0&\le b_2^\pm(\overline X^\e)+\delta(\C X^\e),
\endaligned
\tag1-6
$$
where $\delta(\C X^\e)=\dim\,\ker(H_2(\R X^\e;\R)\to H_2(\C X^\e;\R))$
and $c^+$, $c^-$, $c^0$ are the numbers of the connected
oriented components of $\R X^\e$ with the negative, positive
and zero Euler characteristic respectively
(we refer to \cite{Ar} for the original formulations and arguments,
see also \cite{V4} and subsection 5.1 below).
The Smith theory gives an estimate $\delta(\C X^\e)\le1$ along with
a certain
information about the topology of $\R A$ in the case $\delta(\C X^\e)=1$.

In \cite{Zv}, Zvonilov extended the Arnold inequalities to the
non-singular  curves $\C A$ of odd
degrees applying a version of these inequalities to the curve obtained
from $\C A$ by adding a line;
this involved real curves with nodal singularities.
The case of arbitrary plane real nodal curves was considered
by Viro \cite{V1}.
In this case the numbers of components, $c^\pm$, $c^0$,
in the Arnold inequalities must be replaced
by the inertia indices, $\s_\pm(\psi_{\ee})$, and the nullity,
$\s_0(\psi_{\ee})$, of the complex intersection form, $\psi_{\ee}$,
on $H_2(\R X^\e)$; an estimate for $\delta(\C X^\e)$ is also required.
Viro described this form in combinatorial terms
and gave an estimate for $\delta(\C X^\e)$, which was later
improved by Kharlamov and Viro.
The current version of the Arnold--Viro inequalities for nodal curves
together with a version of such inequalities for nodal surfaces
is formulated in the survey of Kharlamov  \cite{Kh1}
(the proofs, based on, or inspired by the ideas in the notes \cite{V3},
are reproduced in \cite{F1}).

\subheading{1.3. The results: generalized Arnold-Viro inequalities}
Allowing $\C X$ to have more complicated singularities,
one has to deal with tree problems:
to evaluate $b_d^\e(\overline X)$, to give an estimate for
$\delta(\C X^\e)$, and, the most important, to give a suitable
combinatorial description of the complex intersection form.
Theorem 8.1.1 allows to express $b_d^\e(\overline X)$ in terms of
$b_d^\e(\C X)$, provided $\C X$ is a $\Q$-homology manifold.
Under a weaker assumption that $\overline X$ is a $\Q$-homology manifold
near $\R X$, Theorem 8.2.1 yields a formula for $b_d^\e(\overline X)$
if $\C X$ has only ICIS (isolated complete intersection singularities),
see \S3.
The estimates for $\delta(\C X^\e)$ (under various assumptions on $\C X$)
are included in Appendix 1.
An example of arrangements of hyperplanes considered in \S6 shows
that this estimate may admit, however, some improvement.
And, finally, the integration formulae (2-3), (2-4), give a required
combinatorial method of calculation of the form $\psi_{\ee}$,
reducing the problem to evaluation of certain local invariant
of singularities (the canonical quadratic form).
To give a description of these local forms for arbitrary dimension of
$\C X$, seems to be not an easy problem; I present in \S5 some methods
of calculation for the surface singularities.
In particular, I justify such a method, announced in \cite{F3},
for the singularities, which appear by taking suspension over
the real curve singularities.

As it often happens, the most general formulation of the result
is not so convenient in applications as its special versions.
Accordingly, we formulate below two such versions of
the generalized Arnold-Viro inequalities,
and refer to \S3 and \S4 for more general formulations.

Assume first that $\C A\subset\Cp{d}$, $d=2n$, is a real hypersurface of
even degree whose singular locus $\Sing(\C A)$ contains
only isolated singularities. Assume also that the double covering
$X\to\Cp{d}$ branched along $\C A$ is a $\Q$-homology manifold.
Denote by $\mu^\pm_x$ and $\mu^0_x$ the $\pm1$-inertia indices and the
nullity of the Milnor form of a singularity at $x\in \Sing(\C A)=\Sing X$,
and let $\mu^\pm$ and $\mu^0$ denote the sum of
$\mu^\pm_x$ and $\mu^0_x$, respectively, taken
for all $x\in \Sing(\C A)$.
Put $\delta(\C X)=\dim\ker\roman{in}$, where
$\roman{in}\: H_d(\R X;\Q)\to H_d(\C X;\Q)$
is the inclusion homomorphism.

Denote by $\C A^\tau$ a real non-singular perturbation
of $\C A$ and by $X^\tau\to\Cp{d}$ the double covering branched
along $\C A^\tau$.
Then the complex intersection forms, $\psi_\pm$, arising
on $H_d(\R X^\pm)$,  satisfy the following inequalities
$$
\aligned
\s_{-\ee}(\psi_\pm)&\le\frac12(b_d^{-\ee}(X^\tau)-\kapnew
-\frac12\mu^{-\ee})+\min(0,\delta(\C X^\pm)-\s_0(\psi_\pm))\\
\s_{\ee}(\psi_\pm)&\le\frac12(b_d^{-\ee}(X^\tau)+\chi(\R X^\pm)+\kapnew
-\mu^{\ee})-1+\min(0,\delta(\C X^\pm)-\s_0(\psi_\pm))\\
\endaligned
\tag1-7
$$
where $\ee=(-1)^n$ and
$\kapnew\in\{0,1\}$ is the $\roman{mod}\,\,2$ residue of $n$.
These inequalities look certainly not complete unless we provide
an estimate for $\delta(\C X)$.
The following such an estimate is given in this paper
$$\delta(\C X)\le(b_d(\C A;\Z/2)-\nu)+(n-1)$$
where $\nu$ is the rank of the inclusion homomorphism
$H_d(\C A;\Z/2)\to H_d(\Cp{d};\Z/2)$.

In the following theorem we restrict ourselves with the case $d=2$, but
allow more general singularities of a reduced real curve
$\C A\subset\Cp2$, requiring only
that $\overline X^+$ (or $\overline X^-$)
is a $\Q$-homology manifold
in a neighborhood of $\R X^+$ (respectively, $\R X^-$); in particular,
we impose no conditions on the imaginary singularities of $\C A$.
Let $\frak p=\frac12(\mu^++\mu^0)$, that is the sum of the genera,
$\frak p_x=\frac12(\mu_x^++\mu^0_x)$, of all the singular points $x\in\C A$.
Put furthermore $\b=\frac12(\sum_{x}\mu^0_x)$,
where the sum is taken for all $x\in\Sing(\C A)-\Sing(\R A)$.
Denote by $\Sing_0(\C A)$ the set of the {\it essential}
singularities of $\C A$, that is,
$\Sing_0(\C A)=\{x\in\Sing(\C A)\,|\,\mu_x=0\}$,
and let $\C A''=\C A-\Sing_0(\C A)$, $\C A'=\C A''\cup\R A$.

\tm{1.3.1. Theorem}
Assume that $\overline X^\pm$ is a $\Q$-homology manifold in a neighborhood
of $\R X^\pm$. Then
$$
\align
\s_+(\psi_{\pm})&\le\frac12(k-1)(k-2)-\frak p+
\min(\frac12(r-\nu),\frac12\til b_0(\C A''),\til b_0(\C A'),\b)\\
\s_+(\psi_{\pm})+\s_0(\psi_{\pm})&\le\frac12(k-1)(k-2)-\frak p+\b+(r-\nu)\\
\s_-(\psi_{\pm})&\le\frac32k(k-1)+\frac12\chi(\R X^\pm)-\frac12\mu^-+
\min(\frac12(r-\nu),\frac12\til b_0(\C A''),\til b_0(\C A'),\b)\\
\s_-(\psi_{\pm})+\s_0(\psi_{\pm})&\le\frac32k(k-1)+\frac12\chi(\R X^\pm)-
\frac12\mu^-+(r-\nu)\\
\endalign
$$
where $r$ is the number of irreducible components of $\C A$.
\endtm

\tm{1.3.2. Corollary}
If $\C X^\pm$ has no essential singularities (i.e., is
a $\Q$-homology manifold), then
$$
\align
\s_+(\psi_{\pm})&\le\frac12(k-1)(k-2)-\frac12\mu^++\min(0,(r-\nu)
-\s_0(\psi_{\pm}))\\
&-\frac12\mu^++(r-\nu)\\
\s_-(\psi_{\pm})&\le\frac32k(k-1)+\frac12\chi(\R X^\pm)-\frac12\mu^-
+\min(0,(r-\nu)-\s_0(\psi_{\pm}))\\
&-\frac12\mu^-+(r-\nu)\\
\endalign
$$
\endtm

\subheading{1.4. Conventions}
We denote by $b_k(X)$ the k-th Betti number of $X$
and by $b_k(X;\Z/2)$ its $\Z/2$-Betti-number,
that is the rank of $H_k(X;\Z/2)$.
Recall that $\Q$-homology manifolds have all the usual
homology properties of manifolds (with the coefficients group $\Q$).
We denote by $\la\G,\D\ra_X$ the intersection index of oriented $d$-cycles,
$\G$, $\D$, in $X$, provided $X$ is
an oriented compact $2d$-dimensional $\Q$-homology manifold,
(or, at least, such a manifold in a neighborhood of $\G\cap\D$).
We denote by $b_{d}^\pm(X)$ the inertia indices
of the intersection form in $X$, when it is well-defined
(i.e., for even $d$ and $X$ being a $\Q$-homology manifold as above),
by $b_{d}^0(X)$ the nullity of this form
(as $X$ may have a non-empty boundary).

The prepositions $\C$, $\R$ and a bar (e.g., $\C X$, $\R X$,
$\overline X$) are used for the complex point sets, the real parts
and the quotients associated to  real algebraic varieties, as well
as to $\conj$-invariant subsets of such  varieties. To simplify
the notation, we identify $\R X$ with its image in $\overline X$.
$\Sing(X)$ denotes the singular locus of a complex algebraic
variety $X$; moreover, for a subset $A\subset X$, we put
$\Sing(A)=\Sing(X)\cap A$, for instance, $\Sing(\R X)=\Sing(\C
X)\cap\R X$. Whenever the construction uses the metric in $\C X$
or $\R X$, we assume that it comes from the $\conj$-invariant
Fubbini-Study metric in $\Cp{N}\supset\C X$.

Note that the standard construction of a semi-algebraic Whitney
stratification of a real algebraic variety $\C X$, is readily
$\conj$-symmetric, which yields a suitable stratification of
$\overline X$ and provide us with a semi-algebraic stratification
of $\R X$, which is obtained by taking intersection of $\R X$ with
the strata of $\C X$. Taking the connected components of these
intersections, we obtain a refinement, $\frak S$, of the above
stratification, that is used in \S2.

Recall that the Euler characteristic with compact support, which
is denoted by $\chi_c$, is additive (see \cite{GM2}) and can be
used as a measure to integrate the appropriate functions, say,
semi-algebraic functions on an algebraic variety (see \cite{V2}).
In this paper, we integrate functions which take constant values
on the strata of $\frak S$. Indeed, such an integration look more
natural and obvious in the PL-category. Accordingly, we formulate
and prove some results in 2.6 for polyhedra, keeping in mind that
algebraic sets can be triangulated (see \cite{Jo} for the modern
proof and the further references on existence of a triangulation).

\subheading{1.5. Acknowledgements}
I would like to thank O. Ya. Viro, who kindly supplied me with the
hand-written notes \cite{V3} sketching the proof of the formulae in
\cite{V1}. I thank also A.~Degtyarev for numerous useful discussions.

\heading
\S2.
Complex intersection of real cycles in a real algebraic variety
\endheading
\subheading{2.1. Inessential singularities and their canonical
bilinear forms}
Let us call a point of a complex $d$-dimensional algebraic variety
an {\it inessential singularity}, or, briefly,
{a \QR-point}, if its  link is a $\Q$-homology $(2d-1)$-sphere
(that is a $\Q$-homology manifold having rational homologies of
$S^{2d-1}$).
We call a complex algebraic variety \QR-variety if it has only inessential
singular points.
A \QR-curve is obviously topologically non-singular, which easily implies
that for any \QR-variety $X$, its {\it topological singularity}
(that is the set of points, $x\in X$, whose links is not homeomorphic
to $S^{2d-1}$) has codimension $\ge2$.

For a point, $x\in \R X$,
of a real algebraic $d$-dimensional variety, we can choose a small
regular compact $\conj$-invariant neighborhood $\C U_x\subset \C X$
(for instance, $\e$-neighborhood in $\CP{N}\supset \C X$, $0<\e<\!<1$),
put $\C M_x=\d\,\C U_x$,
and call $x$ {\it a \QQ-point} if $\overline M_x$ is a rational homology
$(2d-1)$-sphere.
Note that a real \QR-point is a \QQ-point provided $d$ is even.
A real variety $\C X$ of even dimension will be called
{\it a \QQ-variety} if
all its real points are \QQ-points.
It is not difficult to check that the topological singularity of $\R X$
for a \QQ-variety $\C X$ has codimension $\ge2$.

Given a Whitney stratified pseudo-manifold, $Z$, embedded smoothly
(with respect to each stratum) in a smooth manifold, $Y$, we can
define {\it a vector field tangent to $Z$} as a vector field in
$Y$ defined along $Z$, whose restrictions to the strata of $Z$ are
tangent to these strata.

Note that
for any $x\in\R X$, there exists a tangent to
$\R X$ vector field, $\xi$, defined along $\R M_x$, which is {\it transverse
to $\R M_x$ and outward-directed}. Such $\xi$ can be constructed by
a stratified controlled lift \cite{GM1, p.42} of the vector field
$\frac{\d}{\d r}$ in $\R$ with respect to the distance function
$r\:\R X\to \R$, $r(y)=dist(x,y)$.
The vector field $i\xi$ (where $i=\sqrt{-1}$) is tangent
to $\C M_x$ being normal to $\R X$ (and in particular to $\R M_x$).

Let $\gamma$ and $\delta$ be $(d-1)$-cycles representing some
homology classes
$[\gamma],[\delta]\in H_{d-1}(\R M_x)$, and $\delta_{i\xi}$ is the cycle
in $\R M_x$
obtained by a small shift of $\delta$ in the direction of $i\xi$
(formally speaking, ``a small shift'' is carried by the
local flow of a vector field $i\xi$; the problem of existence of
the flow of stratified controlled vector fields is analyzed in \cite{Ma}).

Then $\delta_{i\xi}\cap \R M_x=\oo$ and we can define a bilinear
form $\l_x\:H_{d-1}(\R M_x)\times H_{d-1}(\R M_x)\to\Q$, which
assigns to the pair $([\gamma],[\delta])$ the linking number
between $\gamma$ and $\delta_{i\xi}$ in $\C M_x$. Note that $\l_x$
is well defined and, in particular, is independent of the choice
of $\xi$, because the set of all vector fields satisfying the
requirement imposed on $\xi$ is convex (with respect to the
obvious linear homotopy). We call $\l_x$ {\it the canonical
quadratic form} associated to a real \QR-point $x$ (Proposition
2.2.3 justifies that $\l_x$ is symmetric). One can treat $\l_x$ as
a local complex intersection form, due to the formula (being just
a re-formulation of the definition of $\l_x$) $$
\l_x([\gamma],[\delta])=\la\Con(\gamma),\Con(\delta_\xi)\ra_{\C
U}\tag2-1 $$ where $\Con(\gamma)$, $\Con(\delta_\xi)$ denote the
cones over the cycles $\gamma$, $\delta_\xi$ in $\C
U_x\cong\Con(\C M_x)$.

Given a \QQ{} point $x\in \R X$ we define similarly a form $\overline\l_x$
on $H_{d-1}(\R M_x)$ assigning the linking number in $\overline M_x$
between $\gamma$ and the image of $\delta_{i\xi}$ in $\overline M_x$.

%%%%%%%%%%%%%%%%        2.2 2.2 2.2 2.2

\heading
2.2. Basic properties of the canonical forms of the inessential
singularities
\endheading
\tm{2.2.1. Proposition}
Assume that $x\in \R X$ is a \QR-point and $d=\dim\C X$ is even. Then
$\overline\l_x=2\l_x$
\endtm

\demo{Proof}
Let  $\gamma$ and $\delta$ be
a pair of $(d-1)$-cycles in $\R M_x$ (considered below also as cycles
in $\C M_x$ and in $\overline M$) with the coefficients in $\Q$.
Then $\gamma=\d\,\overline{\s}$ for some $d$-chain
$\overline{\s}$ in $\overline M_x$, and for the pull back, $\s$,
of $\overline{\s}$ in $\C M_x$ we have obviously $\d\,\s=2\gamma$.
Let $\delta_{i\xi}$ and $\delta_{-i\xi}$ denote the cycles in $\C M_x$
obtained from $\delta$ by shifts in the direction of $i\xi$ and $-i\xi$.
Then
$\la\s,\delta_{i\xi}\ra_{\C M_x}=\la\s,\delta_{-i\xi}\ra_{\C M_x}$,
because $\conj$ preserves the orientation of $M_x$ for even $d$ and
permutes $\delta_{i\xi}$ and $\delta_{-i\xi}$, and thus
$$
2\l_x([\gamma],[\delta])=
\la{\s},{\delta}_{i\xi}\ra_{\C M_x}=
\la\overline{\s},\overline{\delta}_{i\xi}\ra_{\overline M_x}=
\overline\l_x([\gamma],[\delta])
\tag2-2$$
where $\overline{\delta}_{i\xi}$ is the image of ${\delta}_{i\xi}$
in $\overline M_x$.
\qed
\enddemo

For the next property note that the product of \QR-varieties
is again a \QR-variety. Furthermore,
if $\C X$ and $\C Y$ are {\it real} varieties and $x\in \R X$, $y\in \R Y$,
then a homeomorphism between
the real link, $\R M_{z}$, at $z=(x,y)$ and
the join $\R M_x*\R M_y$ of the real links at $x$, $y$ yields
a canonical isomorphism
$H_{p+q-1}(\R M_{z})\cong H_{p-1}(\R M_x)\otimes H_{q-1}(\R M_y)$,
where $p=\dim\R X$, $q=\dim\R Y$.

\tm{2.2.2. Proposition}
Assume that $\C X$ and $\C Y$ are real \QR-varieties and
$z\in \R X\times \R Y$.
Then the canonical form $\l_{z}$ is isomorphic to
$(-1)^{pq}\l_x\otimes\l_y$, where $\l_x$, $\l_y$ are the canonical
forms at the points $x\in\R X$ and $y\in\R Y$, $z=(x,y)$.
\endtm

Combining Propositions 2.2.1 and 2.2.2, we obtain
(for even $p$ and $q$)
$\overline\l_z\cong\frac12\overline\l_x\otimes\overline\l_y$.
Another corollary is that for a non-singular point $x\in \R X$,
$\l_{z}$ differs from $\l_y$ only by the sign.
In particular, if
 $x\in \R X$ belongs to a stratum
$S\subset \R X$ of dimension $p$ and codimension $q$, then
$\l_x\cong(-1)^{pq+\binom{p+1}2}\l_S$, where $\l_S$ denotes the
canonical form of the singularity at $x$ in a real normal slice to
$S$ in $\R X$. By  definition, such a slice is cut on $\R X$ by a
real $(N-p)$-dimensional plane in $\Rp{N}\supset \R X$, which is
transversal to $S$. Since the strata $S$ of our stratification,
$\frak S$, are connected, the form $\l_S$ is independent of the
choice of $x\in S$, due to the local triviality of $\C X$ along
the Whitney strata.

\demo{Proof of Proposition 2.2.2}
Choose $(p-1)$-cycles, $\gamma_1$, $\delta_1$, in $\R M_x$
and $(q-1)$-cycles,
$\gamma_2$, $\delta_2$, in $\R M_y$ and denote by $\delta_1'$,
$\delta_2'$ the cycles obtained from $\delta_1$, $\delta_2$ by small
shifts along the canonical framings of $\R M_x$ and $\R M_y$.
Then $\delta_1'*\delta_2'$ is obtained by a similar shift of
$\delta_1*\delta_2$ in $\C M_z\cong\C M_x*\C M_y$.
Using (2-1), we have
$$
\align
\l_{z}([\gamma_1*\gamma_2],[\delta_1*\delta_2])&=
\la\Con(\gamma_1*\gamma_2), \Con(\delta_1'*\delta_2')
\ra_{{}_{\C X\times \C Y}}\\
&=\la\Con(\gamma_1)\times \Con(\gamma_2),\Con(\delta_1')\times
\Con(\delta_2')\ra_{{}_{\C X\times \C Y}}\\
&=(-1)^{pq}\la\Con(\gamma_1),\Con(\delta_1')\ra_{{}_{\C X}}
\la\Con(\gamma_2), \Con(\delta_2')\ra_{{}_{\C Y}}\\
&=(-1)^{pq}\l_x([\gamma_1],[\delta_1])\l_y([\gamma_2],[\delta_2])
\endalign
$$
\qed
\enddemo

\tm{2.2.3. Proposition}
The form $\l_x$ as well as
$\overline\l_x$ is symmetric, whenever it is well defined.
\endtm

\demo{Proof}
First, note that
$\lk(\gamma,\delta_{i\xi})=(-1)^d\lk(\gamma,\delta_{-i\xi})$,
where $\lk$ stands for the linking number in $\C M_x$.
This is because
$\conj$ changes the orientation for odd $d$, while preserving it for even
$d$, and keeps $\gamma$ fixed, interchanging
$\delta_{i\xi}$ and $\delta_{-i\xi}$.
Furthermore, we have obviously
$\lk(\gamma_{i\xi},\delta)=\lk(\gamma,\delta_{-i\xi})$ and thus
$$
\align
\l_x([\gamma],[\delta])&=\lk(\gamma,\delta_{i\xi})=
(-1)^d\lk(\gamma,\delta_{-i\xi})=
(-1)^d\lk(\gamma_{i\xi},\delta)\\
&=\lk(\delta,\gamma_{i\xi})=\l_x([\delta],[\gamma])
\endalign
$$
In the case of the forms $\overline\l_x$, the arguments are the same,
except that $d$ may not be odd.
\qed
\enddemo

Finally, we consider the natural generators of $H_d(\R X)$ for a real
\QR-variety $\C X$ of any dimension $d$, or \QQ-variety of even dimension,
and give
an integration formula for complex intersection of real cycles in $\C X$.

By a {\it component of $\R X$} we mean the closure of a connected
component of $\R X-\Sing_{\roman{top}}(\R X)$, where
$\Sing_{\roman{top}}(\R X)$ is the topological singularity of $\R X$.
As was mentioned in 2.1, the codimension of $\Sing(\R X)$ is $\ge2$
for both \QR- and \QQ-varieties, so the components of $\R X$
can be viewed as $(\Z/2)$-cycles generating $H_d(\R X;\Z/2)$.
 Similarly, the orientable (outside $\Sing_{\roman{top}}(\R X)$)
components of $\R X$, after we fix an orientation, represent the
generators of $H_d(\R X)$.

\comment
\tm{2.2.4. Lemma}
Assume that $\C X$ is a real \QR-variety.
Then $\Sing_{\top}(\C X)$  has
complex codimension $\ge2$ in $\C X$.
If $\C X$ is a \QQ-variety, then $\Sing_{\top}(\R X)$ has
real codimension $\ge2$ in $\R X$.
\endtm

\demo{Proof}
Assume that $S\subset\Sing(\C X)$ is a codimension 1 stratum in
a \QR-variety $\C X$.
A normal slice to $S$ in $\C X$ at a point $x\in S$ is a germ of a
complex curve, which has a \QR-singularity  at $x$, since the link of $x$
is $\C X$ is an iterated suspension over its link in the normal slice.
 Since the complex link of a \QR-singularity for a curve can be only
$S^1$ (and thus the real link must be $S^0$),
the links $\C M_x$ and $\R M_x$ of $x$ in $\C X$ and $\R X$ must
be also spheres, as the iterated suspension over spheres.
In the case of \QQ-variety, there is in addition a possibility that
the normal-slice-link consists of a pair
of circles permuted by $\conj$.
% (nothing else is possible, since
%the $(-1)$-eigenspace of $\conj_*$ acting in the first homologies
%of this link must be obviously one-dimensional).
It can be easily checked that in the latter case the quotient link
$\overline M_x$ is also a sphere.
\enddemo
\endcomment

Given a pair of oriented components $\G$, $\D$ of $\R X$
and $x\in\G\cap\D$, we put $\l_x(\G,\D)=\l_x([\gamma],[\delta])$,
where $\gamma$, $\delta$ are the cycles on $\R M_x$ cut by $\G$ and $\D$.

\tm{2.2.4. Theorem}
Assume that $\G$ and $\D$ are oriented components of $\R X$ and
$\G\cap\D$ contains only \QR-points. Then
$$
\la\G,\D\ra_{{}_{\C X}}=\int_{\G\cap\D}\l_x(\G,\D)\,d\chi(x)
\tag2-3$$
If $d$ is even and $\R X$ contains only \QQ-points, then
$$
\la\G,\D\ra_{{}_{\overline X}}=\int_{\G\cap\D}\overline\l_x(\G,\D)\,d\chi(x)
\tag2-4$$
\endtm

\tm{2.2.5. Corollary}
If $d$ is even and $\G\cap\D$ contains only \QR-points, then
$$\la\G,\D\ra_{\overline X}=2\la\G,\D\ra_{\C X}$$.
\endtm

\rk{Remark}
Note that $\la\G,\D\ra_{\C X}$ must vanish for any $\G$, $\D$
if the dimension $d$ is odd, however, $\l_x$ may be non-trivial
(see Example 2.4.3).
\endrk

%%%%%%%%%%%%%%%%%%%%%%%%%%%%%%  2.3 2.3 2.3 2.3

\subheading{2.3. Proof of Theorem 2.2.4}
To evaluate $\la\G,\D\ra_{\C X}$ we follow the same approach
as in non-singular case, i.e., shift $\D$ by the flow of $i\eta$,
where $\eta$ is an appropriate vector field,
tangent to the strata of $\R X$.
We suppose that $\eta$ has
a finite set, $\Sigma\subset\R X$, of singularities (zeros),
so that $\G\cap\D_{i\eta}=\Sigma$, where $\D_{i\eta}$ is
a cycle obtained from $\D$ by a shift.
Then $\la\G,\D\ra_{\C X}$ is the sum of the local intersection indices,
$\la\G,\D_{i\eta}\ra_x$, at $x\in\Sigma$.
For any $x\in\Sigma$, we express $\la\G,\D_{i\eta}\ra_x$
in terms of the forms $\l_x$, provided
the flow of $\eta$ is positive (expanding)
in the direction normal to the stratum, $S$, containing $x$.

To give a precise formulation of our assumption about $\eta$,
we recall that
the local triviality of Whitney stratified spaces
along a stratum, $S$, (cf. \cite{GM1, p. 37}), implies that
there exists a chart $\ch_x\:U\to \R^p\times \R^{N-p}$ in a neighborhood
$U\subset\Rp{N}$ of $x\in S\subset\R X$, mapping the points of $\R X$
to $\R^p\times N$ and, in particular, the points of the
stratum $S$ to $\R^p\times\{0\}$.
Here $N$ is a Whitney stratified space called
{\it the normal slice} of $S$ in $\R X$.
We call $\eta$ {\it a stratified expanding vector field with respect to
a stratification}, $\frak S$, of $\R X$, if
for any $x\in\Sigma$ there exists a
``product chart'' as above, in which $\eta$ splits into a direct sum,
$\eta(x,y)=\eta_S(x)+\eta_N(y)$,
where $\eta_S$ is a vector field in $\R^p$ (the components of $\eta$
along $S$) and $\eta_N$ is a field in $\R^{N-p}$, which is required
to have a positive flow.
Denote by $\ind_x(\eta,\frak S)$ the index of $\eta_S$ at $x$
(if $S=\{x\}$ is a $0$-dimensional stratum, then
$\ind_x(\eta,\frak S)=1$).
A singularity of $\eta_S$ will be called {\it elementary} if it is
a standard
non-degenerated singularity in some chart around $x$.
A standard singularity of a vector field in $\R^p$
is by definition represented either by
``the identity'' vector field, $\xi(x)=x$, or, by the
``mirror reflection'' field,
$\xi(x_1,x_2,\dots,x_p)=(-x_1,x_2,\dots,x_p)$.
If for a vector field, $\eta$, on $\R X$,
its restrictions, $\eta_S$, to the strata, $S\in\frak S$
have only elementary singularities, then $\eta$ will be called an
{\it elementary vector field}.

Theorem 2.2.4 follows from the following lemmas.

\tm{2.3.1. Lemma}
Given a real variety $\C X$ endowed with a Whitney stratification $\frak S$,
there exists an elementary stratified expanding vector field, $\eta$,
with respect  to $\frak S$.
\endtm

\tm{2.3.2. Lemma}
Assume that $\eta$ is like in lemma 2.3.1
and $S\subset \R X$ is a stratum of $\frak S$.
Then
$$
\sum_{x\in\Sigma\cap S}\ind_x(\eta,\frak S)=\chi_c(S)
\tag2-5$$
\endtm

\tm{2.3.3. Lemma}
Assume that $\G$ and $\D$ are like in Theorem 2.2.4, whereas
$\eta$ and $S$ are like in Lemma 2.3.1.
Then for any $x\in\Sigma\cap S$
$$
\la\G,\D_{i\eta}\ra_x
=\l_x(\Gamma,\Delta)\ind_x(\eta,\frak S),
\tag2-6$$
\endtm

Using (2-5)-(2-6), we obtain the formula (2-3) as follows
$$
\align
\la\G,\D\ra_{\C X}&=
\sum_{S\subset \G\cap\D}\sum_{x\in S\cap\Sigma} \la\G,\D_{i\eta}\ra_x
=\sum_{S\subset\G\cap\D}\sum_{x\in S\cap\Sigma}\l_x(\Gamma,\Delta)
\ind_x(\eta,\frak S)\\
&=\sum_{S\subset\G\cap\D}\l_S(\G,\D)\chi_c(S)
=\int_{\G\cap\D}\l_x(\G,\D)\,d\chi(x)
\endalign
$$
where we put
$\l_S(\G,\D)=\l_x(\G,\D)$ for $x\in S$, which makes sense, because
$\l_x(\G,\D)$ is independent of $x\in S$ due to
the local triviality of $\C X$ along $S$ and connectedness of $S$.
The proof of (2-4) is analogous.
\qed

\demo{Proof of Lemma 2.3.1}
Denote by $Z_n$ the union of the strata of $\frak S$ of dimension
$(\le n)$.
We construct inductively a smooth vector field $\eta_n$ on
a neighborhood, $V_n\supset Z_n$ in $\Rp{N}\supset\R X$, so that $\eta_n$
coincides with $\eta_{n-1}$
on a neighborhood $V_{n-1}'\subset V_{n-1}\cap V_n$ of $Z_{n-1}$ and
satisfies the properties of a stratified expanding
elementary vector field with respect  to
$\frak S$
 (which include that $\eta_n$ is tangent to $\R X$ along $\R X\cap V_n$,
and have only elementary singularities, which form a set
$\Sigma_n\subset Z_n$).

If $x\in Z_0$, then we let $\eta_0(x)=0$ and define $\eta_0$
around $x$ as the stratified controlled lift of the vector field
$r^2\frac{\d}{\d r}$ on $\R$, where $r\:\R X\to \R$ is
the distance from $x$.
Given $\eta_{n-1}$, it is not difficult to
extend it to $n$-strata, possibly, varying $\eta_{n-1}$
outside some neighborhood of $Z_{n-1}$.
For a generic such an extension, the singularities at
$x\in\Sigma_n\setminus Z_{n-1}$ will be, obviously, non-degenerated.
Furthermore, using an isotopy having support in a small neighborhood of $x$,
it is not difficult to reduce a non-degenerated singularity
to one of the two standard patterns making it elementary.

Finally, using a product chart,
$\ch_x\:U\to \R^n\times\Rp{N-n}$, around $x\in Z_n-Z_{n-1}$,
like above, we can locally extend our vector field, $\eta_S$,
as a direct sum, $\eta_S+\eta_N$, to $U$. The field
$\eta_N$ around the origin in $\Rp{N-n}$ is constructed similarly
to $\eta_0$, using a stratified controlled lift of $r^2\frac{\d}{\d r}$,
with $r$ being the distance function from the stratum $S$
(i.e., from $\R^n$) in the chart $\ch_x$.
In particular, $\eta_N$ has positive flow and is tangent to
$\R X$ at the points of $\R X\cap U$.
Patching together such local extensions
via a partition of unity,
we construct a field $\eta_n$ defined, as is required,  in a
neighborhood, $V_n\supset Z_n$.
Reducing the size of the domain $V_n$, if it is needed,
we can make $\eta_n$ have no zeros in $V_n-Z_n$.
\qed
\enddemo

\demo{Proof of Lemma 2.3.2}
Let $F_A$ denote the local flow of the vector field $\eta$ restricted
to $A\subset \R X$.
The proof consists in applying
the Lefschetz fixed point formula, in the form
\cite{GM2}, to $F_S$, for a stratum, $S$, of $\frak S$, which gives
$$
\chi_c(S)=\sum_{x\in\Sigma\cap S}\Lef_x(F_S)
\tag2-7$$
where $\Lef_x(F_S)$ is the local Lefschetz number of $F_S$ at $x$,
which coincides with $\ind_x(\eta,\frak S)$.
Since $S$ is non-compact, we cannot directly apply
the result of \cite{GM2}, and obtain (2-7) as the difference
of the Lefschetz formulae applied to $F_{\Cl S}$ and $F_{\d S}$,
where $\d S=\Cl S-S$
$$
\aligned
\chi(\Cl S)&=\sum_{x\in\Sigma\cap\Cl S}\Lef_x(F_{\Cl S})\\
\chi(\d S)&=\sum_{x\in\Sigma\cap\d S}\Lef_x(F_{\d S})
\endaligned
\tag2-8$$
Recall that
the Lefschetz formula in \cite{GM2} requires
{\it weak hyperbolicity} of a mapping
(this property generalizes the Morse non-degeneracy
condition and means, informally speaking,
that a mapping can be representing near a fixed point
as an expansion in one direction and a contraction in the complementary one),
which is satisfied at least in the case of primitive stratified expanding
vector fields
with respect to $\frak S$.
To obtain (2--7), it is only left to notice that
$\Lef_x(F_{\Cl S})=\Lef_x(F_{\d S})=\ind_x(-\eta|_T)$, where $-\eta|_T$ is
the restriction of $-\eta$ to the stratum,
$T\subset\d S$, of $\frak S$ containing $x$.
These equalities follow from that the local Lefschetz number
for germ of a mapping, as was defined in \cite{GM2}, does not
change if we take its direct product with
a germ of a contraction.
\qed
\enddemo

\demo{Proof of Lemma 2.3.3}
Consider $x\in\Sigma\cap S$, where $S$ is a stratum of $\frak S$.
If $\ind_x(\eta,\frak S)=1$ and thus $\eta_S$ is ``the identity''
vector field, in some chart around $x$, then we have
$\la\G,\D_{i\eta}\ra_x=\lk_x(\gamma,\delta_{i\eta})=\l_x(\gamma,\delta)$,
so (2--6) holds. Here, like in the previous subsection,
$\gamma$ and $\delta$ are the cycles cut in $\C M_x$ by $\G$ and $\D$,
and $\lk_x$ is the linking number in $\C M_x$.
If $\ind_x(\eta,\frak S)=-1$, then we should replace
``the identity'' vector field by
``the mirror reflection''. But such a change effects on
$\lk_x(\gamma,\delta_{i\eta})$ and $\ind_x(\eta,\frak S)$ as
multiplication by $-1$, so (2-6) still holds, since
$\l_x(\gamma,\delta)$ is preserved, as it is independent of $\eta$.
\qed
\enddemo

%%%%%%%%%%%%%%%%%%%%%%%%%%%%%%%%%%%%%%%%%%%%%

\subheading{2.4. \sQ{} singularities of hypersurfaces and their
canonical forms}
Consider an analytic reduced (i.e., not containing multiple factors)
singularity $f\:(\C^{d},0)\to(\C,0)$,
defining a germ of a hypersurface $\{f=0\}=(\C A,0)\subset(\C^d,0)$.
If $0\in\C A$ is a \QR-point, then $f$ will be called {\it \QR-singularity}.
Similarly, by a {\it \QQ-singularity} we mean
the complexification of a real analytic germ, $f$, such that
$0\in\C A$ is a \QQ-point.
If  {\it the suspension}, $f^S\:(\C^{d+1},0)\to(\C,0)$,
$f^S(x_1,\dots,x_{d+1})=f(x_1,\dots,x_d)-x_{d+1}^2$, over $f$,
is \QR-singularity, (\QQ-singularity),
then we call $f$ a \sQ{\it-singularity}
(respectively, \sQQ{\it-singularity}).

Note that
the projection $\pi\:(\C X,0)\to(\C^d,0)$
forgetting the last coordinate of the level set
$\{f^S=0\}=(\C X,0)\subset(\C^{d+1},0)$
is a double covering branched along $(\C A,0)$ and
the singular locus, $\Sing(\C X)=\Sing(\C A)$, has codimension $\ge2$.
Denote by $\C B\subset\C^d$ a compact $\e$-ball ($0\!<\!\e\!<\!\!<\!1$)
around $0$, let $\C S=\d(\C B)$ and put
$\R B_\pm=\{x\in\R B\,|\,\pm f(x)\ge0\}$,
$\R S_\pm=\R S\cap \R B_\pm$.
Denote by $V_i$, $i=1,\dots,s$, the closures of the connected components
of $\R S\setminus\R A$ and
call $V_i$ {\it local partition  regions of $\R A$ at $0$}, putting
$\sign(V_i)\in\{+,-\}$ for the sign of $f$ inside $V_i$.
Put  $\C M=\pi^{-1}(\C S)$,
and orient $\R M=\pi^{-1}(\R S_+)$ (in the complement of $\Sing(\R X)$),
so that the restriction of the projection $\R M\to\R S$ to
$\R M\cap(\R^d\times\R_+)$
preserves the orientation, whereas its restriction to
$\R M\cap(\R^d\times\R_-)$ reverses
($\R S$ is oriented here as a boundary of $\R B\subset\R^n$).
Similarly we orient $\pi^{-1}(\R S_-)\subset\R^d\times i\R$,
making the restriction of $\pi$ preserve the orientation on
$\pi^{-1}(\R S_-)\cap(\R^d\times i\R_+)$,
and reverse on $\pi^{-1}(\R S_-)\cap(\R^d\times i\R_-)$.
With the inherited orientation,
$\B_i=\pi^{-1}(V_i)$, $i=1,\dots, s$, can be viewed as
oriented cycles in $\C M=\pi^{-1}(\C S)$
(note that $\codim\Sing(\B_i)\ge2$) and we denote by $[\B_i]$ their
fundamental classes.

Given a \sQ-singularity $f$, we define
a $\Q$-valued form  $\l^S$ on $\Cal H=
H^0(\R S\setminus\R A)\cong H_{d-1}(\R S,\R A\cap\R S)$
 using the following version of the construction in 2.1.
Denote by $\xi$ a smooth vector field defined in a neighborhood of $\R S$
in $\R^d$, which is transverse to $\R S$ and
outward-looking along it, being also
tangent to $\R A$ at the points $x\in\R A$
(to construct $\xi$ we use the stratified lifting theorem, like in 2.1).
Denote by $\xi^S$ the vector field in $\R^{d+1}$
tangent to $\pi^{-1}(\R^d)\subset\C X$, obtained by lifting of $\xi$.
Note that the vector field $i\xi^S$ is
tangent to $\C M$  and normal
(in the standard metric of $\C^{d+1}$)
to the set $R_M=p^{-1}(\R S)\subset \C M$
at points $x\in R_M$.

Denote by $v_i$, $i=1\dots,s$, the generators of $\Cal H$ represented by
the {\it characteristic cochains} of $V_i$
(equal to $1$ on $V_i\setminus\R A$ and to $0$ on the rest
of $\R S\setminus\R A$), and put
$$
\l^S(v_i,v_j)=\lk(\B_i,\B_j')
$$
where $\lk$ is the linking number in $\C M$ and
$\B_j'$ is obtained from $\B_j$ by a small shift in the direction
of $i\xi^S$ (i.e., by the corresponding flow).

Note that the opposite choice of
the orientation of $\R^d$ (and, thus, of $\R S$),
changes the orientation of $\B_i$, but does not change the form $\l^S$.
However, changing the sign of $f$, we interchange the summands,
$\Cal H_\pm=H^0(\R S_\pm\setminus\R A)
\cong H_{d-1}(\R S_\pm,\R S_\pm\cap \R A)$,
in the splitting $\Cal H=\Cal H_+\oplus \Cal H_-$.

Denote by $\l^{S}_\pm$ the restrictions of $\l^S$ to $\Cal H_\pm$.
The construction of 2.1 applied to a cone-like neighborhood
$\C U=\pi^{-1}(B)$
of $0\in \C X$ defines the canonical form on $H_{d-1}(\R M)$ and it is
not difficult to see that
$\l^{S+}$ is its pull back via the following product map
$$
H^0(\R S_+\setminus\R A)\cong H^0(\R S_+\setminus\Sing(\R A))@>\pi^*>>
H^0(\R M\setminus\Sing(\R M))\cong H_{d-1}(\R M,\Sing(\R M))
\cong H_{d-1}(\R M)
$$
where the last isomorphism
is due to that $\Sing(\R M)$ has codimension $\ge2$.

Similarly, we define a $\Q$-valued form $\bar{\l}^S$ on $\Cal H_+$,
provided $d$ is even and $f$ is a \sQQ-singularity.
The arguments of Proposition 2.2.1 show that $\bar{\l}^S=2\l^S|_{\Cal H_+}$
if $f$ is a \sQR-singularity and $d$ is even.

\rk{Remark} By the Edmonds theorem, the fixed point
set of a smooth involution on a Spin manifold gets certain
{\it semi-orientation} (that is a pair of the opposite orientations)
provided the involution preserves the orientation and the Spin structure.
If $f$ is an isolated singularity,
then the Milnor fiber, $\C U_t=(f^S)^{-1}(t)$, $t\in \R$, of $f^S$
is {\it Spin} and for $d\ge2$ is simply connected, thus, if $d$ is even, then
$\conj$ preserves both the orientation and
the Spin structure in $\C U_t$ and thus in $\C M_t=\d(\C U_t)$.
This endows $\R M\cong \R M_t$ with a semi-orientation.
A slightly modified version of this construction can be applied to
non-isolated singularities as well, and it is not difficult
to show that the orientation of $\R M$, that we constructed above,
coincides with a $\Spin$ semi-orientation.
\endrk

Note that although the restrictions $\l^{S}_\pm$ of $\l^S$
are quadratic forms (isomorphic to
the canonical quadratic forms of the singularities $f\pm x_{d+1}^2$),
the bilinear form $\l^S$
itself is not symmetric.

\tm{2.4.1. Proposition}
Assume that $\sign(V_i)\ne\sign(V_j)$. Then
$\l^S(v_i,v_j)=-\l^S(v_j,v_i)$.
\endtm

The proof is analogous to that of Proposition 2.2.3.
\qed

This property of $\l^S$ look more natural after changing the basis
in $\Cal H_{\C}=\Cal H\otimes\C$. Namely, we put
$\til v_k=iv_k$ (here $i=\sqrt{-1}$), if $\sign(V_k)=-$, and
$\til v_k=v_k$ if $\sign(V_k)=+$, $k=1,\dots,l$.
The bilinear extension, $\l_\C^S$, of $\l^S$ to $\Cal H_{\C}$
is defined then by a self-adjoint matrix and
the following ``product formula'' holds.

Assume that $h\:(\C^{p+q},0)\to(\C,0)$ is a product singularity,
$h(x,y)=f(x)g(y)$, where $f\:(\C^{p},0)\to(\C,0)$ and
$g\:(\C^{q},0)\to(\C,0)$ are \QS-singularities. Then
$h$ is also a \QS-singularity, since
the germ $X_h=\{h^S=0\}$ is obviously isomorphic
to the quotient $(X_f\times X_g)/\theta$, where
$X_f$ and $X_g$ are  the germs defined by the equations $f^S=0$ and $g^S=0$
and $\theta$ is
the direct product of the deck transformations of the branched coverings
$(X_f,0)\to(\C^p,0)$ and $(X_g,0)\to(\C^q,0)$.
Denote by $V_1',\dots,V_s'$ and $V_1'',\dots,V_t''$ the local partition
regions for $f$ and $g$ respectively.
Then such regions for $h$ correspond to $V_{kl}=V_k'*V_l''$ under
the natural homeomorphism $\R S_h\cong \R S_f*\R S_g$
(the objects like $\R S$, $\Cal H$, $\Cal H_\C$ and $\l^S_\C$
are marked by the subscripts $f$, $g$, $h$ if they are
associated to the corresponding singularities).
Let $v_k'\in\Cal H_f$, $v_l''\in\Cal H_g$, $v_{kl}\in\Cal H_h$
and $\til v_k'\in\Cal H_{f\C}$, $\til v_l''\in\Cal H_{g\C}$,
$\til v_{kl}\in\Cal H_{h\C}$ denote the corresponding bases.
Consider the isomorphism
$\Cal H_{f\C}\otimes_{\C}\Cal H_{g\C}\cong \Cal H_{h\C}$
which sends
$\til v_k'\otimes \til v_l''$ to $\til v_{kl}$.
The arguments analogous to that of Propositions 2.2.1--2.2.2
prove the following relation.

\tm{2.4.2. Proposition}
$\l^S_{h\C}\cong\frac12(-1)^{pq}\l^S_{f\C}\otimes \l^S_{g\C}$
\qed\endtm
{\bf 2.4.3. Example.}
The forms $\l^S$ and $\l^S_\C$ for the identity function
$f\:(\C,0)\to (\C,0)$, $f(x)=x$, are defined by the matrices
$\pmatrix
1&-1\\
1&\phantom{-}1\\
\endpmatrix$
and
$\pmatrix
1&-i\\
i&-1\\
\endpmatrix$
in the bases $\{v_k\}$ and $\{\tilde v_k\}$ ($k=1,2$)
respectively.

Using Proposition 2.4.2 we can determine the form $\l^S_\C$ for
$f\:(\C^d,0)\to (\C,0)$, $f(x_1,\dots,x_d)=x_1\dots x_d$.
Namely, let us mark the local partition regions of $f$ by
vectors $a=(a_1,\dots,a_d)$, $a_k\in\{+1,-1\}$ belonging to these
regions. Let $v_a\in\Cal H$, and $\til v_a\in\Cal H_\C$ denote
the basis elements representing the region containing $a$. Then
$$
\aligned
\l^S_\C(\til v_a,\til v_b)
&=(-1)^{\frac{d(d-1)}2}2^{1-d}\sign(b)\,i^{d(a,b)}
=(-1)^{\frac{d(d-1)}2}2^{1-d}\sign(a)\,(-i)^{d(a,b)}\\
\l^S(v_a,v_b)&=(-1)^{\frac{d(d-1)}2}2^{1-d}\sign(b)\,
i^{d(a,b)-a_--b_-}\\
\endaligned
\tag2-9
$$
where $\sign(a)=a_1\dots a_d$, and $a_-$, $b_-$ denote the number of
negative coordinates, whereas $d(a,b)$ the number of
distinct coordinates in vectors $a$ and $b$.
%Namely, $\l^S(v_i,v_j)=2^{-d+1}(-1)^{\binom{d}2+c}$,
%where $c$ is the codimension of $V_i\cap V_j$ in $\R S$.

Finally, we introduce a relative version of the form $\l^S$,
involving a ($\Z/2$)-cycle $\Omega\subset\R^d$ of codimension 1, which
contains $0$ and have smooth simplices transversally intersecting $\R S$.
Transversality yields a $(d-1)$-cycle,
$\Omega_S=\Omega\cap \R S$, in $\R S$. Its complement
$\R S-\Omega_S$ splits, by the Alexander duality, in two regions
distinguished by $(\roman{mod}\ 2)$-linking number with $\Omega_S$.
Assuming that $\Omega_S$ does not intersect the interiors of
$V_i$, $V_j$, we put
$\l^S(v_i,v_j|\Omega)=\l^S(v_i,v_j)$, if the both $V_i$, $V_j$ lie
in the closure of one of the above two regions, and
$\l^S(v_i,v_j|\Omega)=-\l^S(v_i,v_j)$ if they lie in
the closures of distinct regions.

Similarly, we define $\bar{\l}^S(v_i,v_j)$ and $\bar{\l}^S(v_i,v_j|\Omega)$,
for $v_i,v_j\in\Cal H_+$, if $f$ is a \sQQ-singularity.

%%%%%%%%%%%%%%%%%%%%%%%%%%%%%%%%%%%%%%%%%%

\subheading{2.5. The partition components of real hypersurfaces}
Consider a non-singular real  variety $\C P$, of dimension $d$
and a real line bundle $\ell\: \C\Cal L\to \C P$, that is
a line bundle supplied with an anti-linear involution,
$\conj_{\Cal L}\: \C\Cal L \to \C\Cal L$, commuting with $\ell$
and the complex conjugation in $\C P$.
Let $\ell_{{}_\R}\:\R\Cal L\to \R P$ denote the real part of $\ell$,
that is its restriction
to the real parts of  $\C\Cal L$ and $\C P$.

The complex conjugation, $\conj_\Cal L{\otimes}\conj_\Cal L$,
makes the square $\C\Cal L{\otimes}\C\Cal L$ a real bundle, whose
real part is trivialized by choosing  the direction of
the ``positive'' ray
$\{x\otimes x\| x\in\R\Cal L\}$ in each fiber.
Thus, for a real (i.e., conjugation-equivariant) section
$f: \C P\to \C\Cal L{\otimes}\C\Cal L$,
the sign of $f(x)$ is well defined at the real points $x\in \R P$.
We put $\R P_\pm=\R P_\pm(f)= \{ x\in \R P\,:\,\pm f\ge 0\}$
and assume in what follows that the zero locus
$\C A\subset \C P$ of $f$ is a reduced hypersurface.

Denote by $\C X\subset \C\Cal L$ the pull back of
$f(\C P)\subset \C\Cal L^{\otimes2}$ via
the mapping $\C\Cal L\to \C\Cal L^{\otimes2}$,
$x\mapsto x\otimes x$.
The restriction $\pi=\ell|_{\C X}\:\C X\to \C P$
is obviously a double covering branched along $\C A$.

%%%%%%%%%%%%%%%%%%%%%%%%%%%

Denote by $W_j$, $j=1,\dots,m$, the closures of the connected
components of $\R P-\R A$
and call $W_j$ the {\it partition component} of $\R A$,
putting $\sign(W_i)$ for the sign of $f$ inside $W_i$.
Consider the class
$\om=w_1(\R P)+w_1(\ell_\R)\in H^1(\R P;\Z/2)$
and note that the restriction $\om|_{W_j}$ vanishes if and only if
$\G_i=\pi^{-1}(W_i)\subset\C X$ is orientable (in the complement
of $\Sing(\C X)$).
This is because $w_1(\R\Cal L)=\ell^*(\om)$ and
the normal bundle to $\G_i-\Sing(\G_i)$ is trivial
(if $\sign(W_i)=+$, then we consider
the normal bundle in $\R\Cal L\supset\G_i$, otherwise we consider it
in $i\R\Cal L\supset\G_i$).

%%%%%%%%%%%%%%%%%%%%%%

Let $W_1,\dots,W_l$, $l\le m$, be the components $W_i$ with
the restriction $\om|_{W_j}=0$,
and put $\R P^\circ=\bigcup_{i=1}^l\Int W_i$.
Realize the homology class dual to $\om$ by a
$\Z/2$-cycle in $\R P$ with smooth simplices
and denote by $\Omega$ the union of the simplices.
Then there exists an orientation of
$\R\Cal L-\ell_\R^{-1}(\Omega)$ which cannot be extended through
the ``walls'' of $\ell_\R^{-1}(\Omega)$, and such an orientation is unique
up to the natural action of $H^0(\R P;\Z/2)$.
It is not difficult to choose $\Omega$ having support in
$\R P-\R P^\circ$.
If we fix such an orientation of $\R\Cal L-\ell_\R^{-1}(\Omega)$,
defined by $\Omega$,
and consider the orientation of $i(\R\Cal L-\ell_\R^{-1}(\Omega))$
induced from it by the homeomorphism $\R\Cal L\to i\R\Cal L$,
$x\mapsto ix$,
then $\G_j$, $j=1,\dots,l$, become oriented cycles.

Given $x\in\R P$, we mark with a subscript $x$
the objects $\R S$, $\l^S$, $\Cal H$, etc., introduced in subsection 2.4,
which are associated to the germ of $f$ at $x$.
Consider {\it the natural basis}, $w_i\in H^0(\R P-\R A)$, \
$i=1,\dots,m$, represented by the characteristic cochains
of $W_i\setminus\R A$ and
let  $w_i(x)\in\Cal H_x$ denote the image of $w_i$ under the inclusion
homomorphism $H^0(\R P-\R A)\to H^0(\R S_x\setminus\R A)$.
We call $f$ a
{\it \sQ-section} and $\C A$ {\it a \sQ-hypersurface}
(respectively, {\it a \sQQ-section} and {\it a \sQQ-hypersurface})
if $\Sing(\C A)$ contains only \sQ-singularities
(\sQQ-singularities); this is obviously
equivalent to that $\C X$ is a \QR-variety (\QQ-variety).
If $\C A$ is a \sQ-hypersurface, then
we define a bilinear $\Q$-valued
{\it partition form} $\phi$ on $H^0(\R P^{\circ})$ putting
$\phi(w_i,w_j|\Omega)=\la\G_i,\G_j\ra_{\C X}$, $1\le i,j\le l$
(here we keep the same notation, $w_i$, for the restriction of
$w_i$ to $H^0(\R P^{\circ})$).
Note that $\phi$ is well defined, in spite of
the ambiguity in the choice of the orientation
of  $\R\Cal L-\ell_\R^{-1}(\Omega)$.

\tm{2.5.1. Theorem}
Assume that $\C A\subset\C P$ is a \sQ-hypersurface. Then
$$
\aligned
\phi(w_i,w_j&|\Omega)=
\int_{W_i\cap W_j}\l^S_x(w_i(x),w_j(x)\|\Omega)\,d\chi(x)
\phantom{AAAAAAA} \text{if }\ 1\le i,j\le l, i\ne j\\
\phi(w_i,w_i&|\Omega)=
\int_{W_i\cap \Sing(\R A)}\l^S_x(w_i(x),w_i(x)\|\Omega)\,d\chi(x)\\
&+(-1)^{\frac{d(d-1)}{2}}(2\chi_c(W_i\setminus\R A)
+\chi_c((\R A\cap W_i)\setminus\Sing(\R A))
\phantom{A}\text{if }\ 1\le i\le l.
\endaligned\tag2-10
$$
\endtm

\demo{Proof}
It is not difficult to see that
the orientations of the components $\G_i$ induce the same
orientations of the cycles $\gamma_i=\R M_x\cap\G_i$ as
the components $\B_j$ in the previous section.

If $\sign(W_i)=\sign(W_j)$ are positive or negative,
then the formulae of Theorem 2.5.1 follow from the formulae of
Theorem 2.2.4.
In the case of the opposite signs the proof is analogous.
\qed
\enddemo

 Note that for even $d$ the form $\phi$ is symmetric and splits into
a direct sum, $\phi=\phi_+\oplus\phi_-$, in
$H^0(\Rp{\circ})=H^0(\Rp{\circ}_+)\oplus H^0(\Rp{\circ}_-)$,
where $\Rp{\circ}_\pm=\Rp{\circ}\cap\Rp{}_\pm$.
If $f$ is a \sQQ-section, then
we define a form $\bar{\phi}\:H^0(\Rp{\circ}_+)\to\Q$
putting $\bar{\phi}(w_i,w_j)=\la\G_i,\G_j\ra_{\overline X}$.
Theorem 2.2.5 and the formula (2-4) then imply the following formula
analogous to (2-10).
$$
\aligned
\bar{\phi}(w_i,w_j&|\Omega)=
\int_{W_i\cap W_j}\bar{\l}^S_x(w_i(x),w_j(x)\|\Omega)\,d\chi(x)
\phantom{AAAAAAA}
\text{if } 1\le i,j\le l, i\ne j\\
\bar{\phi}(w_i,w_i&|\Omega)=
\int_{W_i\cap \Sing(\R A)}\bar{\l}^S_x(w_i(x),w_i(x)\|\Omega)\,d\chi(x)\\
&+(-1)^{\frac{d(d-1)}{2}}(4\chi_c(W_i\setminus\R A)
+2\chi_c((\R A\cap W_i)\setminus\Sing(\R A))
\phantom{A}\text{if }\ 1\le i\le l.
\endaligned\tag2-11
$$
Theorem 2.2.1 implies also that
$\bar{\phi}(w_i,w_j)=2\phi_+(w_i,w_j)$, provided
 $\C A$ is a \sQR-hypersurface and $d$ is even.

%%%%%%%%%%%%%%%%%%%%%%%%%%%%%%%%%%%%%%%%%

\subheading{2.6. Residue form}
In this section we show how integration along the odd-dimensional strata
in the formulae (2-10), (2-10),
can be reduced to integration along their boundary.
Assume that $Q$ is a compact polyhedron of dimension $d$.
Let us call a function $f\:Q\to\R$ {\it constructible}
if it is constant on the open simplices of some triangulation,
$\Cal T$, of $Q$.
For such a function, we may consider its restriction to the link,
$\Lk_x(Q)$, of $x\in Q$ and define
$$
\hat f(x)=\int_{\Lk_x(Q)}f(y)\,d\chi(y)
$$
The latter definition makes sense if the link $\Lk_x(Q)$ is taken
with respect to a sufficiently fine triangulation,
say, the barycentric subdivision of
any refinement of $\Cal T$, containing $x$ as a vertex
(here, as above, $f$ is constant of the simplices of $\Cal T$).
Alternatively, we may assume that $\Lk_x(Q)$ is defined as the
infinitesimal link of $x$ (the direct limit of the usual
links of $x$ with respect to all triangulation containing $x$
as a vertex, or, equivalently, the set of germs of $PL$-rays with
the origin at $x$), and define the restriction $f|_{\Lk_x(Q)}$
in the obvious way.

\tm{2.6.1. Lemma}
For any constructible function $f$ on a compact polyhedron $Q$
$$
\int_{Q}\hat f(x)\,d\chi(x)=0
$$
\endtm

\demo{Proof}
This identity can be easily checked if $f$ is a characteristic
function of a closed simplex of $\Cal T$.
In general, $f$ is a linear combination of such functions and
the formula of the lemma follows from  additivity of the integral.
\qed\enddemo

\tm{2.6.2. Corollary}
Assume that the dimension, $d$, of $Q$
is odd and $Q_{d-1}$ denotes the union of the
$k$-simplices, $k\le d-1$, of
$\Cal T$. Then
$$
\int_{Q} f(x)\,d\chi(x)=\int_{Q_{d-1}}
(f(x)-\frac12\hat f(x))\,d\chi(x)
$$
\endtm

\demo{Proof}
This follows from Lemma 2.6.1 and from that
$f(x)-\frac12\hat f(x)=0$ inside $d$-simplices of $\Cal T$,
for odd $d$.
\qed\enddemo

Let $Q_{\roman{sing}}$ denote {\it the topological singularity} of $Q$,
that is the set of points $x\in Q$ whose link, $\Lk_x(Q)$,
is not homeomorphic to $(d-1)$-sphere, where $d=\dim Q$.

\tm{2.6.3. Corollary}
Assume that $Q$ is a compact polyhedron of dimension $d$
and $Q_{\roman{sing}}\subset Q'\subset Q$ for a sub-polyhedron $Q'$. Then
$$
\align
\int_{Q'}\chi(\Lk_x(Q))\,d\chi(x)&=0, \phantom{AAAAAAAAA}
\text{if $d$ is even}\\
\int_{Q'}\chi(\Lk_x(Q))\,d\chi(x)&=-2\chi_c(Q-Q'), \phantom{AA}
\text{if $d$ is odd}
\endalign
$$
\endtm

\demo{Proof}
It follows from Lemma 2.6.1 applied to the constant function $f=1$,
since it gives $\hat f(x)=\chi(\Lk_x(Q))=1-(-1)^d$ for
$x\in Q-Q_{\roman{sing}}$.
\qed\enddemo

\comment
\rk{Remarks}
\roster
\item
Lemma 2.6.1 and its corollaries can be applied to Whitney stratified
spaces and the functions constant of the strata, since such spaces are
known to admit a triangulation.
\item
In Corollary 2.6.3 we may equally assume that $Q'$ is just
any union of simplices of any triangulation of $Q$,
rather then a sub-polyhedron.
\endroster
\endrk
\endcomment

Consider a real reduced hyper-surface $\C A$ in a real nonsingular
$d$-dimensional variety $\R P$ and a partition component $W_i$
defined as in the section 2.5.
For $x\in\R P$, let
$\chi_x(W_i)=\chi(W_i\cap\R S_x)$
where $\R S_x$ is  an $\e$-sphere, $0<\e<\!<1$, around $x$
in $\Rp{d}$.

Applying Corollary 2.6.3 for $Q=W_i$ and $Q'=\R A\cap W_i$, we obtain

\tm{2.6.4. Corollary}
Assume that $\C A$ is like above. Then
$$
\align
&\int_{W_i\cap\Sing(\R A)}\chi_x(W_i)\,d\chi(x)=
-\chi_c((\R A\cap W_i)\setminus\Sing(\R A))
\phantom{AAAAAAAAAA}\text{if $d$ is even}\\
&\int_{W_i\cap\Sing(\R A)}\chi_x(W_i)\,d\chi(x)=
-2\chi_c(W_i\setminus\R A)-\chi_c((\R A\cap W_i)\setminus\Sing(\R A))
\phantom{A}
\text{if $d$ is odd}\\
\endalign
$$
\qed\endtm

Using Corollary 2.6.2 and 2.6.4 we can rewrite the formulae (2-10)
as follows (the ``hat'' over $\l^S_x$ has below the same meaning as
it has over $f$ in Lemma 2.6.1).

\tm{2.6.5. Corollary}
Assume that $\C A$ is like in Theorem 2.5.1. Then
$$
\aligned
\phi(w_i,w_j|\Omega)=\int_{(W_i\cap W_j)_{k-1}}
(\l^S_x(w_i(x),w_j(x)&|\Omega)-
\frac12\hat{\l}^S_x(w_i(x),w_j(x)|\Omega))
\,d\chi(x)\\
\text{ if }  &i\ne j \text{ and } k=\dim(W_i\cap W_j) \text{ is odd},\\
\phi(w_i,w_i|\Omega)=\int_{W_i\cap\Sing(\R A)}
(\l_x^S(w_i(x),w_i(x)&|\Omega)-
(-1)^{\frac{d(d-1)}2}\chi_x(W_i))\,d\chi(x)\\
+(-1)^{\frac{d(d-1)}{2}}&(1+(-1)^d)\chi_c(W_i\setminus\R A)
\phantom{A}\text{if}\  1\le i\le l.
\endaligned
$$
Here $(W_i\cap W_j)_{k-1}$ is constituted by the points of $W_i\cap W_j$
which belong to the strata of dimension $\le k-1$.
\endtm

Note that Corollary 2.6.3 can be applied to $W_i\cap W_j$, since its
topological singularity is contained in $(W_i\cap W_j)_{k-1}$.
\qed

The above formula becomes more simple if we assume that $x$ is an
{\it isolated} \sQ-singularity. In this case, the term
$\chi_x(W_i)$ can be understood as the value
$\chi_x(w_i(x),w_i(x))$ of the bilinear form, $\chi_x\:\Cal
H_x\times \Cal H_x\to\Z$, defined in the basis
$v_1,\dots,v_s\in\Cal H_x$ as follows $$ \align
\chi_x(v_i,v_i)&=\chi(V_i)\\
\chi_x(v_i,v_j)&=\sign(v_i)\frac12\chi(V_i\cap V_j),\ \text{if}\
i\ne j\\
\endalign
$$ where $\sign(v_i)=+1$ if $V_i$ is a positive local partition
component and $\sign(v_i)=-1$ if negative. The term under the
integral in Corollary 2.6.5 then becomes $\q(w_i(x),w_j(x))$,
where $\q\:\Cal H_x\times \Cal H_x\to\Q$ is the bilinear form
$$\q=\l_x^S-(-1)^{\frac{d(d-1)}2}\chi_x$$ which will be called
{\it the residue form}. We define also the relative form,
$\q(\centerdot,\centerdot|\Omega)$, of $\q$ using the same
convention as in the subsection 2.4 for the form $\l$.

The formulae (2-10) can be rewritten as follows

\tm{2.6.6. Corollary}
Assume that $\C A\subset\C P$ is a hypersurface, like in Theorem 2.5.1,
which have only isolated \sQR-singularities. Consider a pair,
$W_i$, $W_j$ of partition components.
Then
$$
\aligned
\phi(w_i,w_j|\Omega)=&
\sum_{x\in W_i\cap \Sing(\R A)}\q(w_i(x),w_j(x)\|\Omega)\\
&+(-1)^{\frac{d(d-1)}{2}}(1+(-1)^d)
\chi_c(W_i\setminus\R A)\,\delta_{ij}
\phantom{AAA}\text{if }\ 1\le i\le l.
\endaligned
\tag2-12 $$ where $\delta_{ij}$ is the Kronecker symbol.
\endtm

Similar arguments give an analogous formula for the form
$\bar\phi$ in the case of \sQQ-singularities. Namely, we consider
the form $\bar q=\bar\l_x^S-2(-1)^n\chi_x|_{\Cal H_x^+}$ on $\Cal
H_x^+$ together with its relative form defined as usual. For
\sQR-singularities we have $\bar q=2\q_+$, where $\q_{\pm}=
\q|{\Cal H_\pm}$. The formulae (2-11) can be stated then as
follows.

\tm{2.6.7. Corollary} Assume that $\C A\subset\C P$ is a
\sQQ-hypersurface, with isolated singularities, $d=2n$ and
$\sign(W_i)=\sign(W_j)=+$. Then $$ \aligned
\bar\phi(w_i,w_j|\Omega)=& \sum_{x\in W_i\cap \Sing(\R A)} \bar
q(w_i(x),w_j(x)\|\Omega)\\ &+2(-1)^n \chi_c(W_i\setminus\R
A)\,\delta_{ij} \phantom{AAA}\text{if }\ 1\le i\le l.
\endaligned
\tag2-13
$$
\endtm

In case of \sQ-hypersurface, for even
$d$ and $\sign(W_i)\ne\sign(W_j)$, we have obviously
$\phi(w_i,w_j)=0$, so, (2-12) implies that
$$
0=\phi(w_i,w_j)=
\sum_{x\in\Sing(\R A)\cap W_i\cap W_j}\q(w_i(x),w_j(x)|\Omega)
$$
It is not difficult to derive from the latter that
all the terms of the above sum must vanish. This implies furthermore that
$\q(v_i(x),v_j(x))=0$, if $\sign(V_i)\ne\sign(V_j)$, for
any isolated singularity $x\in\R A$, and thus
$\l^S_x(v_i,v_j)=(-1)^\frac{d(d-1)}2\chi_x(v_i(x),v_j(x))$.
We summarize it as follows.

\tm{2.6.8. Proposition}
Assume that $f\:(\C^{2n+1},0)\to(\C,0)$ is an isolated \sQ-singularity
and $\l^S$, $\q$ are the forms associated to it as above.
Then
\roster
\item
$\l^S(v_i,v_j)=\sign(v_i)(-1)^{n}\frac12\chi(V_i\cap V_j)$
if $\sign(V_i)\ne\sign(V_j)$.
\item
the form $\q$ splits into a direct sum $\q=\q_+\oplus \q_-$.
\endroster
\endtm
%%%%%%%%%%%%%%%%%%%%%%%%%%

%       C O M P L E T E     I N T E S E C T I O N S

%%%%%%%%%%%%%%%%%%%%%%%%%%
\heading \S3. Generalized Arnold--Viro inequalities for complete
intersections
\endheading
\subheading{3.1. The results} Given a real algebraic variety, $\C
X$, of dimension $d$, we let $\delta(\C X)=\dim\,\ker(H_d(\R
X;\R)\to H_d(\C X;\R))$, which is obviously equal to
$\dim\,\ker(H_d(\R X;\R)\to H_d(\overline X;\R))$. As was
mentioned in the introduction, the generalized Arnold--Viro
inequalities are the estimates $$ \align \s_\pm(\psi)&\le
b_d^\pm(\overline X)\\ \s_\pm(\psi)+\s_0(\psi)&\le
b_d^\pm(\overline X)+\delta(\C X)\\
\endalign
$$ for the inertia indices of the complex intersection form $\psi$
in $H_d(\R X)$, being expressed in some suitable form. We can
evaluate $b_d^\pm(\overline X)$ in terms of $b_d^\pm(\C X)$, as it
is done in Appendix B (Theorem 8.1.1). For instance, in the case
of real complete intersection \QR-varieties of dimension $d=2n$,
we obtain $$ \aligned b_d^{-\ee}(\overline
X)&=\frac12(b_d^{-\ee}(\C X)-\kapnew)\\ b_d^{\ee}(\overline
X)&=\frac12(b_d^{\ee}(\C X)+\chi(\R X)-\kapnew)
\endaligned\tag3-1
$$ where $\ee=(-1)^{n}$ and $\kapnew=\frac12(1-\ee)$. Furthermore,
for even $n$ we have actually an estimate
$\s_\kapnew(\psi)=\s_+\le b_d^\kapnew(\overline X)-1$. To see it,
note that the hyperplane section class, $H\in H_{d-2}(\C X)$, is
anti-invariant with respect to $\conj$ (i.e., $\conj_*(H)=-H$) and
thus must be orthogonal to the image of $H_*(\R X)$ in the
Lefschetz ring, $H_*(\C X)$. In particular, $h^n$ vanishes on
$H_d(\R X)$, where $h\in H^2(\C X;\Q)$ is dual to $H$. On the
other hand, $h^n$ is $\conj^*$-invariant for even $n$, and thus
descends to a positive-square class in $H^d(\overline X;\Q)$.

These results can be  summarized in the  following theorem.

\tm{3.1.1. Theorem} Assume that $\C X$ is a real complete
intersection \QR-variety of dimension $d=2n$ and
 $\ee=(-1)^{n}$, $\kapnew=\frac12(1-\ee)$.
Then $$ \align \s_{-\ee}(\psi)&\le\frac12(b_d^{-\ee}(\C
X)-\kapnew) +\min(0,\delta(\C X)-\s_0(\psi)) \tag3-2\\
\s_{\ee}(\psi)&\le\frac12(b_d^{\ee}(\C X)+\chi(\R X)+\kapnew)-1+
\min(0,\delta(\C X)-\s_0(\psi))\tag3-3\\
\endalign
$$ \qed
\endtm

Evaluation of $b_d^{-\ee}(\C X)$ for an arbitrary \QR-variety is
beyond the scope of this paper. We only note that if such a
variety, $\C X$, has only isolated complete intersection
singularities (ICIS), then $b_d^\pm(\C X)=b_d^\pm(\C
X^\tau)-\mu^\pm$, where $\C X^\tau$ denote a perturbation of $\C
X$ (that is a non-singular real complete intersection obtained by
a small variation of the equations defining $\C X$), and $\mu^\pm$
are the total Milnor numbers for $\C X$ defined like in 1.3.
Evaluation of $b_d^\pm(\C X^\tau)$, in the case of complete
intersections, is also an easy problem, since the known Chern
classes of  $\C X^\tau$ determine obviously both $\chi(\C X^\tau)$
and $\s(\C X^\tau)$.

Under a weaker assumption that $\C X$ is a \QQ-variety, a result
similar to Theorem 3.1.1 needs less trivial calculations. We
remove from $\C X$ a regular $\conj$-symmetric compact regular
neighborhood, $\C U_0$, of $\Sing_0(\C X)-\Sing_0(\R X)$ (the
purely imaginary essential singularity) let $\C X'=\Cl(\C X-\C
U_0)$, and follow a similar approach applying it to the quotient
$\overline X'$, which is a $\Q$-homology manifold. The result,
which we present here only for the case of ICIS, is as follows.

\tm{3.1.2. Theorem} Assume that $\C X$ is a complete intersection
\QQ-variety of dimension $d=2n$, with only isolated singularities.
Then $$ \aligned
\s_{-\ee}(\overline{\psi})&\le\frac12(b_d^{-\ee}(\C
X^\tau)-\kapnew)-\frak p +\min(\gamma,\b+\delta(\C
X')-\s_0(\psi))\\
\s_{\ee}(\overline{\psi})&\le\frac12(b_d^{-\ee}(\C X^\tau)+\chi(\R
X) +\kapnew-\mu^{\ee})-1+\min(\gamma-\b,\delta(\C
X')-\s_0(\psi))\\
\endaligned\tag3-4
$$ where $\b=b_{d}(\d \overline U_0)$, $\gamma=b_{d+1}(\overline
X)$, $\delta(\C X')=\dim\,\ker(H_d(\R X;\R)\to H_d(\overline
X';\R))$ and $\frak p=\frac12(\mu^{-\ee}+\mu^0)$.
\endtm

To make the inequalities in Theorems 3.1.1--3.1.2 usable, we need
to complete them by estimating $\delta(\C X)$ and $\delta(\C X')$.
Like in the case of the usual Arnold inequalities, such estimates
come from the Smith theory (for the proof see Appendix A).

\tm{3.1.3. Proposition} Assume that $\C X$ is a real projective
algebraic variety of dimension $d$ or a $\conj$-invariant subset
of such a variety. Then $$ \delta(\C X)\le\sum_{k=d+1}^{2d}b_k(\C
X;\Z/2)-b_{d+1}(\overline X) \tag3-5$$ In particular, if $\C X$ is
a real complete intersection having only isolated singularities,
and $d=2n$, then $$\delta(\C X)\le b_{d+1}(\C
X;\Z/2)+n-b_{d+1}(\overline X)\tag3-6$$
Assume that $\C X\to\Cp{d}$, $d=2n$, is a double covering
branched along a real reduced hypersurface $\C A\subset \Cp{d}$.
Then $$ \delta(\C X)\le\sum_{k=d+1}^{2d-1}b_k(\Cp{d},\C
A;\Z/2)+(n-1) -b_{d+1}(\overline X), \tag3-7 $$ If $\C A$ has  in
addition only isolated singularities, then $$\delta(\C
X)\le(b_d(\C A;\Z/2)-\nu_d)+n-b_{d+1}(\overline X)\tag3-8$$
\endtm

\rk{Remarks} \roster
\item
The estimate (3-7) is better by 1 then (3-8) for $d=2$.
\item
If $\C X$ is a \QR-variety and $d=\dim\C X$ is even, then
$b_{d+1}(\overline X)=\frac12 b_{d+1}(\C X)=\frac12 b_{d-1}(\C X)$
(one can show it using that the mixed Hodge structure in $H^*(\C
X)$ is pure, cf. 8.3). If moreover $\C X$ is a complete
intersection, then $b_{d+1}(\overline X)=0$.
\item
The estimates (3-7) and (3-8) still hold if $\Cp{d}$ is replaced
by any complete intersection real non-singular variety of
dimension $d$, and a hypersurface $\C A$ is very ample (viewed as
a divisor).
\endroster
\endrk

\demo{Proof of Theorem 3.1.2} The following Proposition 3.1.4
(proved in subsection 3.3) evaluates $b_d^\pm(\overline X')$ and
$b_d^0(\overline X')$. To estimate $\delta(\C X')$ we use (3-5).
The rest is analogous to the proof of Theorem 3.1.1. \qed\enddemo

\tm{3.1.4. Proposition} Assuming that $\C X$ is like in Theorem
3.1.2, we have $$ \align b_d^0(\overline X')&=\b-\gamma\tag3-9\\
b_d^{-\ee}(\overline X')&=\frac12(b_d^{-\ee}(\C
X^\tau)-\kapnew)-\frak p +\gamma\tag3-10\\ b_d^{\ee}(\overline
X')&=\frac12(b_d^{\ee}(\C X^\tau)+\chi(\R X)-\kapnew
-\mu^{\ee})+\gamma-\b\tag3-11\\
\endalign
$$
\endtm

\subheading{3.2. Two properties of the real ICIS (isolated
complete intersection singularities)} For the proof of Proposition
3.1.4 we need to use two results about real ICIS. The first one is
a version of the Milnor Lemma for the quotients by the complex
conjugation of a real ICIS. Given such a singularity
$f\:(\C^{d+k},0)\to(\C^k,0)$, we put $\C U=\{f=0\}\cap\C
B\subset\C^{d+k}$, where $\C B$ is a compact $\e$-ball around
zero, $0<\e<\!<1$, and consider a small real deformation
$f^\tau\:(\C^{d+k},0)\to(\C^k,0)$, $0<|\tau|<\!<\e$, of $f=f^0$
along with the corresponding deformation of $\C U$, denoted by $\C
U^\tau$. We call $f^\tau$ (along with $\C U^\tau$) {\it a
perturbation}, if $\C U^\tau$ is non-singular. It is well known
(cf. \cite{Lo}, \cite{D}) that $\C U^\tau$ is homotopy equivalent
to a wedge of $d$-spheres, for any deformation, $f^\tau$.

\tm{3.2.1. Lemma} Assume that $\C U^\tau$ is a deformation of a
cone-like compact neighborhood, $\C U$, for a real ICIS of
dimension $d\ge1$. Then the quotient $\overline U^\tau$ is
homotopy equivalent to a wedge of $d$-spheres, provided $\R
U^\tau\ne\oo$. If $\R U^\tau=\oo$, then $\overline U^\tau$ still
has rational homology of a wedge of $d$-spheres.
\endtm

\demo{Proof} Note that for $d=1$ the statement of the lemma is
trivial, and the condition $\R U^\tau\ne\oo$ implies that
$\overline U^\tau$ is simply connected for $d\ge2$. Furthermore,
$\pi_k(\overline U^\tau)=0$ for $2\le k\le d-1$, since a generic
mapping of $S^k$ to $\overline U^\tau$ does not intersect $\R
U^\tau$ and can be lifted to $\C U^\tau$. In addition,
$H_k(\overline U^\tau;F)=0$ for $k\ne d$, where $F$ is a field of
the characteristic $\ne2$,
 for instance $\Q$ or $\Z/p$, for a prime $p\ne2$
(as it is well known that the projection $\C U^\tau\to \overline
U^\tau$ induces an isomorphism between $H^k(\overline U^\tau;F)$
and the $\conj_*$-invariant subspace of $H^k(\C U^\tau;F)$).

It follows also from the Smith sequence for $\conj$ (see 7.1) that
$H_k(\overline U^\tau,\R U^\tau;\Z/2)= H_{k+1}(\overline U^\tau,\R
U^\tau;\Z/2)$ for $k\ge d+1$. Thus, $H_k(\overline
U^\tau;\Z/2)=H_k(\overline U^\tau,\R U^\tau;\Z/2)=0$, for $k\ge
d+1$. Furthermore,  $H_d(\overline U^\tau)$ is torsion free by the
universal coefficients formula, since $H_{d+1}(\overline
U^\tau;\Z/p)=0$ for any prime $p$, so it is left to apply the
Whitehead theorem.
\enddemo

\rk{Remark} One can apply the same arguments in a more general
setting, for instance, for non-isolated real singularities, and
prove that $\overline U^\tau$ and of $\d \overline U\cong\d
\overline U^\tau$ have the same connectedness properties as $\C
U^\tau$ and $\d \C U$ respectively, unless $\conj$ acts freely
(the connectedness properties of $\C U^\tau$ and $\d \C U$ can be
found, e.g., in \cite{Di, p.76}). \comment Note that for $d=2$ the
quotient $\overline U^\tau$ is not only a $\Q$-homology manifold
(as it is for any even $d$), but an honest a oriented 4-manifold,
which can be endowed naturally with a smooth structure.
\endcomment
\endrk

The next result is proven in Appendix B.

\tm{3.2.2. Theorem} Assume that $\C U^\tau$ is a non-singular
perturbation of a real ICIS of dimension, $d=2n\ge2$. Then $$
\align b_d^{-\ee}(\overline U^\tau)+b_d^0(\overline U^\tau)
&=\frac12(b_d^{-\ee}(\C U^\tau)+b_d^0(\C U^\tau))\tag3-12\\
b_d^{\ee}(\overline U^\tau)&=\frac12(b_d^{\ee}(\C U^\tau)
+{\chi}(\R U^\tau)-1) \tag3-13\endalign $$
\endtm

\subheading{3.3. Proof of Proposition 3.1.4} In this subsection we
denote by $\C U\subset \C X$ a regular compact $\conj$-invariant
neighborhood of the whole $\Sing(\C X)$ and by $\C U^\tau\subset
\C X^\tau$ its non-singular perturbation, so that $\C X^\circ=\C
X-\Int(\C U)$ is identified with $\C X^\tau-\Int \C U^\tau$. Note
that $\bar{ X}^\circ$ is obtained from a $\Q$-homology manifold
$\bar X'$ by removing several $\Q$-homology $4$-discs, and thus,
$b_d^\pm(\bar X')=b_d^\pm(\bar{X}^\circ)$ and $b_d^0(\bar
X')=b_d^0(\bar{X}^\circ)$.

Furthermore, it follows from the long exact sequence of $(\bar
X^\circ,\d \bar X^\circ)$ that $$b_d^0(\bar X^\circ)=b_{d-1}(\d
\bar X^\circ)-b_{d-1}(\bar X^\circ)+ b_{d-1}(\bar X^\circ,\d \bar
X^\circ)-\til b_{d-2}(\d\bar X^\circ) $$ since $b_{d-2}(\d \bar
X^\circ)=b_{d-2}(\bar U^\tau)$ vanish for $d>2$ (cf. \cite{Di, p.
76}. This implies (3-9), because, by the duality, $b_{d-1}(\bar
X^\circ)=b_{d+1}(\bar X^\tau,\bar U^\tau)=b_{d+1}(\bar X)$ and
$b_{d-1}(\d \bar X^\circ,\d \bar X^\circ)= b_{d-1}(\bar X)+\til
b_{d-2}(\d\bar X^\circ)$, for $d\ge2$.

The long exact sequence of $(\bar X^\tau,\bar U^\tau)$ yields $$
\align b_d(\bar X^\circ)&= b_d(\bar X^\tau,\bar U^\tau)=b_d(\bar
X^\tau)-b_d(\bar U^\tau)+ b_{d+1}(\bar X^\tau,\bar
U^\tau)-b_{d+1}(\bar X^\tau)\\ &=b_d(\bar X^\tau)-b_d(\bar
U^\tau)+\gamma \tag3-14\endalign $$ Furthermore, $2b_d^{-\ee}(\bar
X^\circ)=b_d(\bar X^\circ)-b_d^0(\bar X^\circ) -\ee\s(\bar
X^\circ)$, where $\s(\bar X^\tau)=\s(\bar X^\circ)+\s(\bar
U^\tau)$, and $$ 2b_d^{-\ee}(\bar X^\circ)=b_d(\bar
X^\tau)-b_d(\bar U^\tau)+ \g-(\b-\g) -\ee(\s(\bar X^\tau)-\s(\bar
U^\tau)) $$

Using (3-12) we obtain $b_d(\bar U^\tau)-\ee\s(\bar
U^\tau)=2b_d^{-\ee}(\bar U^\tau) +b_d^0(\bar U^\tau)=2p-\b$, which
gives $2b_d^{-\ee}(\bar X^\circ)= (b_d(\bar X^\tau)-\ee\s(\bar
X^\tau))+2\gamma-\b-(2p-\b)$ that is $$b_d^{-\ee}(\bar
X^\circ)=b_d^{-\ee}(\bar X^\tau)-p+\gamma \tag3-15$$ The relation
(3-3) applied to $\C X^\tau$ gives (3-10). Subtracting (3-9) and
(3-15) from (3-14), we obtain $$ b_d^{\ee}(\bar
X^\circ)=b_d^{\ee}(\bar X^\tau)-b_d^{\ee}(\bar U^\tau)
+\gamma-\b\tag3-16$$

Finally, (3-13) gives $\chi(\R X^\tau)-\chi(\R X)=\til\chi(\R
U^\tau)= 2b_d^{\ee}(\bar U^\tau)-\mu^\ee$, which together with
(3-16) and (3-3) applied to $\C X^\tau$ implies the last identity
(3-11).

\heading \S4. Generalized Arnold--Viro inequalities for real
algebraic surfaces
\endheading

\subheading{4.1. Generalized Arnold--Viro inequalities for curves}
Let $\rho$ denote the number of real branches at $0$ of the zero
locus, $(\C A,0)=\{f=0\}$, of a real singularity
$f\:(\C^2,0)\to(\C,0)$. Let us call $f$ (as well as $(\C A,0)$)
{\it a dot singularity} if $\rho=0$. A dot singularity can be {\it
positive}, if $f$ is positive around $0$ on $\R^2$, and {\it
negative} otherwise.

Given a real even \sQ-curve, $\C A$, on a non-singular real
surface, $\C P$, we consider the partition form, $\phi$, defined
in 2.5. Recall that $\C A$ is the zero locus of a real section,
$f$, of $\C \Cal L^{\otimes2}$ for some real line bundle,
$\ell\:\C\Cal L\to\C P$, and $f$ being fixed defines a splitting
$\phi=\phi_+\oplus\phi_-$. Recall also that the double covering,
$\pi\:\C X\to\C P$, branched along $\C A$ is the restriction of
$\ell$ to $\C X\subset\C\Cal L$. Denote by $K$ the canonical class
of $\C P$ and by $L$ the divisor class of the line bundle
$\ell\:\C \Cal L\to \C P$, introduced in 2.5.

Assume that $f$ admits a real perturbation, $f^\tau$, so that the
corresponding perturbation,
 $\C A^\tau$, of $\C A$ is non-singular. Let
$\pi^\tau\:\C X^\tau\to\C P$ denote the corresponding perturbation
of $\pi$. For an $\e$-neighborhood, $\C B\subset\C P$, of
$\Sing(\C A)$, we put $\C U=\pi^{-1}(\C B)$, $\C
U^\tau=(\pi^{\tau})^{-1}(\C B)$, assuming that
$0<\!|\tau|\!<\!\!<\!\e<\!\!<\!1$. We consider moreover an
$\e$-neighborhood $\C B_0$ of the essential imaginary singular
locus, $\Sing_0(\C A)\setminus\Sing(\R A)$, and put similarly $\C
U_0=\pi^{-1}(\C B_0)$, $\C U_0^\tau=(\pi^{\tau})^{-1}(\C B_0)$.
Put furthermore $\C P^\circ=\Cl(\C P-\C B)$, $\C A^\circ=\C
A\cap\C P^\circ$, $\C X^\circ=\Cl(\C X-\C U)$ and define similarly
$\C P'$, $\C A'$, $\C X'$, (respectively, $\C P''$, $\C A''$, $\C
X''$) removing from $\C P$, $\C A$, $\C X$ $\e$-neighborhoods of
the essential imaginary singular loci (respectively, neighborhoods
of the whole essential singular loci).

Recall that $\mu^\pm=b_2^\pm(\C U^\tau)$, $\mu^0=b_2^0(\C
U^\tau)$,
 $\frak p=\frac12(\mu^++\mu^0)$.
The imaginary essential singularities of $\C A$ split into
complex-conjugated pairs; denote by $\a^{(0)}_{\Im}$ the number of
these pairs. Consider a very good real resolution $\C
X^{\res}\to\C X$ of $\C X$, denote its exceptional divisor by $\C
E$, and put $\b=b_1(\bar E)$

Denote by $\nu$ ($\nu'$) the rank of the inclusion homomorphism
from $H_2(\C A;\Z/2)$ (respectively, from $H_2(\C A';\Z/2)$) to
$H_2(\C P;\Z/2)$), and let $t_2$ denote the rank of $\Z/2$-torsion
in $H_1(\C P)$.

\tm{4.1.1. Theorem} Assume that $\C A$ is a \sQQ curve. Then $$
\align
%%%%%%%%%%%%%% 4-1      4-1
\s_+(\bar\phi)\le& b_2^+(\C P)+\frac12L.(K+L)-\frak p+t_2\\
&+\frac12\min(\tilde b_0(\C A''),\tilde b_0(\C A)+2\b,b_2(\C
A)-\nu) \tag4-1\\ \s_+(\bar\phi)+\s_0(\bar\phi)\le& b_2^+(\C P)
+\frac12L.(K+L)-\frak p+2t_2+b_1(\C P)+(b_2(\C A')-\nu')+\b\\
&+\max(0,3\a^{(0)}_{\Im}-1)\tag4-2\\ \s_-(\bar\phi)\le& b_2^-(\C
P)+\frac12(L.(K+3L)+\chi(\R X)-\mu^-)+t_2-\b\\
&+\frac12\min(\tilde b_0(\C A''),\tilde b_0(\C A)+2\b,b_2(\C
A)-\nu) \tag4-3\\ \s_-(\bar\phi)+\s_0(\bar\phi)\le & b_2^-(\C
P)+\frac12(L.(K+3L)+\chi(\R X)-\mu^-)+2t_2+b_1(\C P)\\ &+(b_2(\C
A')-\nu')+\max(0,3\a^{(0)}_{\Im}-1)\tag4-4\\
\endalign
$$
\endtm

If $\C A$ is a \QS-curve, then $\a^{(0)}_{\Im}=\b=0$, $\C A=\C
A'=\C A''$. $\frak p=\frac12\mu^+$ and Theorem 4.1.1 can be
applied to estimate $\s_\pm$, $\s_0$ of the both
$\phi_+=\frac12\bar\phi$ and $\phi_-$. This simplifies the
formulae (4-1)---(4-4) as follows

\tm{4.1.2. Corollary} Assume that $\C A$ is a \QS-curve and
$b_1(\C P;\Z/2)=0$, then $$ \align \s_+(\phi_{\e})\le& b_2^+(\C
P)+\frac12L.(K+L)-\frac12\mu^+
+\min(r-\nu-\s_0(\phi_{\e}),0)\tag4-5\\ \s_-(\phi_{\e})\le&
b_2^-(\C P)+\frac12(L.(K+3L)+\chi(\R X)-\mu^-)
+\min(r-\nu-\s_0(\phi_{\e}),0)\tag4-6\\
\endalign
$$ where $r$ is the number of irreducible components of $\C A$ and
$\e\in\{+,-\}$.
\endtm

\demo{Proof of Theorem 4.1.1} Following the scheme of the proof of
Theorem 3.1.2, we need only to interpret $$\align
\s_\pm(\phi_{\e})&\le b_2^\pm(\overline X')\\
\s_\pm(\phi_{\e})+\s_0(\phi_{\e})&\le b_2^\pm(\overline X')
+b_2^0(\overline X') +\delta(\C X')
\endalign$$
modifying the left-hand side in accord with the identities of
Lemma 4.1.3 and the estimates of $\delta(\C X')=\dim\ker(H_2(\R
X)\to H_2(\overline X'))$ in Lemma 4.1.4. \qed\enddemo

%%%%%%%%%%%%%%%%%%%%%%%%%%%%%%%%%%%%%%%%%%%%

\tm{4.1.3. Lemma} Let $\gamma=b_3(\overline X)-b_1(\overline X)$.
Then $$ \align b_2^0(\overline{
X}')=b_2^0(\overline{X}^\circ)=&\b-\g\tag4-7\\
b_2^+(\overline{X'})=b_2^+(\overline{ X}^\circ)=& b_2^+(\C
P)-b_1(\C P)+\frac12L.(K+L)-\frak p+b_3(\overline X)\tag4-8\\
b_2^-(\overline{X'})=b_2^-(\overline{ X}^\circ)=& b_2^-(\C
P)-b_1(\C P)+\frac12(L.(K+3L)+\chi(\R X)-\mu^-)\tag4-9\\
&+b_3(\overline X)-\b\\
\endalign
$$
\endtm

\tm{4.1.4. Lemma} Assume that $\C X\to\C P$ is a morphism of real
surfaces being a double covering  branched along a real reduced
curve $\C A\subset\C P$, where $\C P$ is non-singular. Then $$
\align
%%%%%%%%  4-10   4-10
b_1(\overline X)&=\frac12b_1(\C X)\le b_1(\C P;\Z/2)+ \frac12\til
b_0(A)\tag4-10\\
%%%%%%%%  4-11   4-11
b_3(\overline X)&=\frac12b_3(\C X)
\le
b_1(\C P;\Z/2)+ \frac12\min(\til b_0(\C A''),\til b_0(\C
A)+2\b,b_2(\C A)-\nu) \tag4-11\\ b_3(\overline X)&\le b_1(\C
P;\Z/2)+\til{b}_0(\C A')\\ \delta(\C X')&\le 2b_1(\C
P;\Z/2)-b_1(\overline X)+(b_2(\C A')-\nu')
+\max(0,3\a^{(0)}_{\Im}-1)\tag4-12\\ b_2^0(\bar X)+\delta(\C X')&
\le 2b_1(\C P;\Z/2)-b_3(\overline X)+(b_2(\C A')-\nu')+\b
+\max(0,3\a^{(0)}_{\Im}-1) \tag4-13\\
\endalign
$$
\endtm
The proof of Lemma 4.1.4 is given in Appendix A.

\demo{Proof of Lemma 4.1.3} Like in the proof of Proposition
3.1.4, we have $b_2^\pm(\overline X')=b_2^\pm(\overline{X}^\circ)$
and $b_2^0(\overline X')=b_2^0(\overline{X}^\circ)$. Similarly, we
obtain (4-1) together with the relations $$ \aligned
b_2^+({\overline X}^\circ)=&\frac12(b_2^+(X^\tau)-1)-\frak p+
b_3({\overline X})-b_1({\overline X})\\ b_2^-(\overline{
X}^\circ)=& \frac12(b_2^-(X^\tau)-1+\chi(X_\R)-\mu^-)-\b
+b_3({\overline X})-b_1({\overline X})
\endaligned\tag4-14
$$ where $b_1({\overline X})=\frac12b_1({\C X})$ by (4-10).

Comparing the Riemann-Hurwitz formula for the covering $\C
X^\tau\to \C P$ with the formula for the signature of an
involution applied to the covering transform, we obtain, like in
the case $\C P=\Cp2$ considered in the introduction, $$ \align
\frac12(b_2^+(\C X^\tau)-b_1(\C X^\tau)+1) &=b_2^+(\C P)-b_1(\C
P)+1 -\frac18(2\chi(\C A^\tau)+\la \C A,\C A\ra_{\C P})\tag4-15\\
\frac12(b_2^-(\C X^\tau)-b_1(\C X^\tau)+1) &=b_2^-(\C P)-b_1(\C
P)+1 -\frac18(2\chi(\C A^\tau)-\la \C A,\C A\ra_{\C P})\tag4-16\\
\endalign
$$ where, by the adjunction formula, $$ \align -\frac18(2\chi(\C
A^\tau)+\la \C A,\C A\ra_{\C P}) &=\frac12L(K+L)\\
-\frac18(2\chi(\C A^\tau)-\la \C A,\C A\ra_{\C P})
&=\frac12L.(K+3L)\\
\endalign
$$
%$2\chi_a(L)=-L.(K+L)$ and
%$\chi_a(L)-\frac12\chi(\C A^\tau)=\frac12L(K+3L)$
Combining (4-14) with (4-15) and (4-16), we obtain (4-8) and
(4-9). \qed\enddemo

\subheading{4.2. Generalized Arnold--Viro inequalities for
surfaces} Consider a very good resolution $\res\:\C X^{\res}\to\C
X$, of a \QQ-surface $\C X$ and put $\C U^{\res}=\res^{-1}(\C U)$,
where $\C U$ is a regular neighborhood of $\Sing(\C X)$. Let us
call a point, $x\in\Sing(\C X)$, {\it $\Z/2$-inessential} if its
link is a $\Z/2$-homology sphere, and {\it $\Z/2$-essential}
otherwise. $\Z/2$-essential imaginary singular points in $\C X$
split into conjugated pairs, whose number we denote by
$\a_{\Im}^{(2)}$.
Consider $\C U_{\re}\subset\C U$ consisting of the connected
components of $\C U$ around the real singularities of $\C X$ and
denote by $d^{(2)}_{\re}$ the rank of $\Z/2$-torsion of $H_1(\d\,
\C U_{\re})$. It is not difficult to check that $d^{(2)}_{\re}$
can be equivalently defined as the nullity of the  $(\roman{mod}\
2)$-intersection-form in $\C U^{\res}_{\re}=\res^{-1}(\C
U_{\re})$, or as the rank of $\Discr\otimes\Z/2$, where $\Discr$
is the discriminant group of the lattice $H_2(\C U^{\res}_{\re})$.
Hence, $d_{\re}^{(2)}$ can be easily computed as soon as we know
the resolution graph of $\C X^{\res}$.

\tm{4.2.1. Theorem} Assume that $\C X$ is a real \QQ-surface with
the partition form $\phi$. Then $$ \align
%%%%%%%%%%%   4-17  4-17   4-17
\s_+(\phi)&\le p_g(\C X^{\res})\tag4-17\\
%%%%%%%%%%%   4-18  4-18   4-18
\s_+(\phi)+\s_0(\phi)&\le \chi_a(\C X^{\res})+b_1(\C
X^{\res};\Z/2)+\b+
\max(0,\a_{\Im}^{(2)}-1)+d_{\re}^{(2)}\tag4-18\\
%%%%%%%%%%%   4-19  4-19   4-19
\s_-(\phi)&\le\frac12(b_2^-(\C X^{\res})-1 +\chi(\R X)+\hat\chi(\R
E)) -b_2(\overline E)\tag4-19\\
%%%%%%%%%%%   4-20  4-20   4-20
\s_-(\phi)+\s_0(\phi)&\le \frac12(b_2^-(\C X^{\res}) +1+\chi(\R
X)-\hat\chi(\C E))+b_1(\C X^{\res};\Z/2)\\ &-\frac12b_1(\C
X^{\res})+ \max(0,\a_{\Im}^{(2)}-1)+d_{\re}^{(2)}\tag4-20\\
\endalign
$$ where $\hat\chi(Z)=\chi(Z)-b_0(Z)$ is the reduced Euler
characteristic (in our particular case, $Z=\C E$ or $Z=\R E$).
\endtm

\demo{Proof} We follow again the standard scheme of proving the
Arnold-Viro-type inequalities, using the following lemmas

\tm{4.2.2. Lemma} In the assumptions of Theorem 4.2.1 we have $$
\align b_2^0({\overline X}^\circ)&=\b-\g\tag4-21\\
b_2^+({\overline X}^\circ) &=b_2^+({\overline X}^{\res})=p_g(\C
X^{\res})\tag4-22\\ b_2^-(\overline{X}^\circ) &=b_2^-(\bar
X^{\res})-b_2(\bar U^{\res}) =\frac12(b_2^-(\C X^{\res})-1+\chi(\R
X)+\hat\chi(\R E)) -b_2(\overline E)\tag4-23\\
\endalign
$$
\endtm

\tm{4.2.3. Lemma} For any real surface $\C X$ with normal
singularities we have $$ \align
%%%      4--24  4--24
%%%
\delta(\C X')&\le b_1(\C X^{\res};\Z/2)-\frac12b_1(\C X)+
\max(1,\a_{\Im}^{(2)}\,)+d_{\re}^{(2)}\tag4-24\\
%%%       4--25  4--25
%%%
b_2^0(\overline{X}^\circ)+\delta(\C X')&\le b_1(\C X^{\res};\Z/2)
-\frac12b_1(\C X^{\res})+\b
+\max(1,\a_{\Im}^{(2)}\,)+d_{\re}^{(2)}\tag4-25\\
\endalign
$$
\endtm

So, (4-17) and (4-19) are immediate corollaries of (4-22) and
(4-23), whereas (4-18) follows from (4-17) and (4-25), since
$\chi_a(\C X^{\res})=p_g(\C X^{\res})-\frac12 b_1(\C X^{\res})+1$.
(4-20) follows from (4-19) and (4-25), since $\frac12\hat\chi(\R
E)-b_2(\bar E)+\b=-\frac12\hat\chi(\C E)$. \qed\enddemo

The proof of Lemma 4.2.3 is given in Appendix A.

\demo{Proof of Lemma 4.2.2} The relation (4-21) is a version of
(3-9) (a minor difference in the setting is not essential for the
proof). The following proof of (4-22) and (4-23) is also similar
to the proof of Theorem 3.3.1 (with $\C X^\tau$ being replaced by
$\C X^{\res}$ because the singularities of $\C X$ may be not
ICIS).

Using the duality and the excision theorem, we obtain
$b_2(\overline X^\circ)=b_2(\overline
X^{\res},\overline{U}^{\res})$, which together with the exact
sequence of $(\overline X^{\res},\overline{U}^{\res})$ gives $$
b_2(\overline X^\circ)=b_2(\overline
X^{\res})-b_2(\overline{U}^{\res})+
b_1(\overline{U}^{\res})-b_1(\overline X^{\res}) +b_1(\overline
X^{\res},\overline{U}^{\res})-\til b_0(\overline{U}^{\res}) $$
because the inclusion homomorphism $H_2(\overline U^{\res})\to
H_2(\overline X^{\res})$ is monomorphic (since $\overline
U^{\res}$ has non-degenerated intersection form). Using that
$b_1(\overline{U}^{\res})=b_1(\overline E)=\b$, $b_1(\overline
X^{\res},\overline{U}^{\res})-\til b_0(\overline{U}^{\res})
=b_1(\overline X)$, and $b_1(\overline X^{\res})=b_3(\overline
X^{\res})=b_3(\overline X)$ (note that the relativization
homomorphism $H_3(\overline X^{\res})\to H_3(\overline
X^{\res},\overline{U}^{\res})$ is isomorphism), we obtain $$
\aligned b_2(\overline X^\circ)&=b_2(\overline
X^{\res})-b_2(\overline{U}^{\res}) +\b-\g\\ b_2^+(\overline
X^\circ)+b_2^-(\overline X^\circ)&= b_2(\overline
X^{\res})-b_2(\overline{U}^{\res})\\
\endaligned\tag4-26
$$

Furthermore, $\s(\overline X^{\res})=\s(\overline
X^\circ)+\s(\overline{U}^{\res})= \s(\overline
X^\circ)-b_2(\overline{U}^{\res})$, since
$b_2(\overline{U}^{\res})=b_2^-(\overline{U}^{\res})$. This yields
$$ 2b_2^+(\overline X^\circ)=b_2(\overline
X^\circ)-b_2^0(\overline X^\circ)+ \s(\overline
X^{\res})+b_2(\overline{U}^{\res})= b_2(\overline
X^{\res})+\s(\overline X^{\res})=2b_2^+(\overline X^{\res}) $$ and
(4-22) follows, since $b_2^+(\overline X^{\res})=\frac12(b_2^+(\C
X^{\res})-1)=p_g(\C X^{\res})$ (cf. (1-4)).

The first equality in (4-23) is obtained from (4-22) and (4-26),
whereas the second one uses (1-5) (or (3-3)). \qed\enddemo

\subheading{4.3. Sharpness of the  generalized Arnold-Viro
inequalities} In this subsection we characterize the gaps between
the left hand side and the right hand side in the generalized
Arnold-Viro inequalities. Here I confine myself with the case of
\sQ-curves and assume in addition $b_1(\C P)=0$ and $b_0(\C A)=0$
(this simplify the formulations, but for the arguments it is not
very essential). Note that these conditions imply that $b_1(X)=0$
for $X\to\C P$ being as above. Let us put $$ \align
\D_\e^+&=b_2^+(\C P)+\frac12(L.(K+L)-\mu^+)-\s_+(\phi_\e)\\
\D_\e^-&=b_2^-(\C P)+\frac12(L(K+3L)+\chi(\R
X^\e)-\mu^-)-\s_-(\phi_\e)\\ \D_\e^0&=(r-\nu)-\s_0(\phi_\e)
\endalign
$$ where $\D_\e^0$ can be negative, although $\D_\e^\pm$ and
$\D_\e^\pm+\D_\e^0$ cannot, by Corollary 4.1.2.

Let $g(\C A)$ denote the geometric genus of the curve $\C A$, and
$g_a(\C A)$ the arithmetic genus. Let $n=b_0(\R P_\e-\R A)$ be the
number of the partition components, $W_i$, and $n_\om$ the number
of those, which are not involved in the partition forms
$\phi_\pm$, i.e., for which $\omega|_{W_i}\ne0$. Denote by $2\bra$
the total number of the imaginary branches at the singular points
of $\C A$ (the number of branches at the imaginary singularities
plus the number of the imaginary branches at the real
singularities). Let $\a_{\roman{\Im}}$ denote the number of pairs
of the imaginary singularities of $\C A$ and $\a_+$ the number of
positive dot-singularities

\tm{4.3.1. Proposition} In the above assumptions on the surface
$\C P$ and curve $\C A$, we have $$ \D_\e=b_2(\C P)-\nu+g(\C
A)+\bra-(\a_{\roman{\Im}}+\a_+)-b_1(\Int \R P_\e) +n_\om+b_2(\R
P_\e)$$
\endtm

\tm{4.3.2. Corollary} The generalized Arnold-Viro inequalities are
equalities if $\C P=\Cp2$, $\C A$ splits into rational irreducible
components one of which has odd degree and all the singularities
of $\C A$ are real and have only real branches.
\endtm

\tm{4.3.3. Corollary} Under the assumption of Corollary 4.3.2 the
radical of the partition form $\phi_\e$ has rank $r-1$. In
particular, the matrix of $\phi_\e$ is singular, since $r\ge2$.
\endtm

For a singularity at $x\in\Sing(\C A)$ the Milnor formula
\cite{Mi, Theorem 10.5} gives a relation
$\delta_x=\frac12(\mu_x+r_x-1)$, where $\mu_x$ is the Milnor
number, $r_x$ is the number of branches of $\C A$ at $x$ and
$\delta_x$ is the maximal number of nodes which can appear after a
deformation of this singularity. Denote by $\rho_x$ the number of
{\it real} branches at $x\in\Sing(\R A)$. We have the following
relations for an irreducible curve $\C A$ in a non-singular
surface $\C P$.

$$\align 2\sum_{x\in \Sing(\C A)}\delta_x-\sum_{x\in \Sing(\C
A)}\mu_x &=\sum_{x\in \Sing(\C A)}(r_x-1)\tag4-27\\ \sum_{x\in
\Sing(\C A)}(r_x-1)-\sum_{x\in \Sing(\R A)}(\rho_x-1)
&=2\bra-2\a_{\roman{\Im}}\tag4-28\\ \chi(\R A)&=-\sum_{x\in
\Sing(\R A)}(\rho_x-1)\tag4-29\\ g_a(\C A)+(r-1)&=g(\C
A)+\sum_{x\in \Sing(\C A)}\delta_x\tag4-30\\
\endalign
$$

(4-27) and (4-28) are obvious and (4-29) is a straightforward
estimate of $\chi(\R A)$. (4-30) is the Pl\"ucker formula, see
\cite{Se, p.74} and \cite{Mi, Property 10.5}. In turn, we note
that existence of \AG-morsifications \cite{AC,GZ} for real curve
singularities (cf. \S5) makes the proof of (4-30) elementary,
reducing it to the trivial case of nodal curves, because
\AG-morsifications obviously preserve $r$, $g$ and $\sum_{x\in
\Sing(\C A)}\delta_x$.

\demo{Proof of Proposition 4.3.1} Adding the expressions defining
$\D_\e^+$, $\D_\e^-$ and $\D_\e^0$, we obtain $$ \D=b_2(\C
P)+L.(K+2L)+\frac12\chi(\R X^\e)-\frac12\mu+(r-\nu)-(n-n_\om) $$
The adjunction formula yields $L.(K+2L)=g_a(\C A)-1$, whereas
(4-27)---(4-30) imply that $$ g_a(\C
A)+(r-1)-\frac12\mu+\frac12\chi(\R A)=g(\C
A)+\bra-\a_{\roman{\Im}} $$ Furthermore, obviously
 $\frac12\chi(\R X^\e)=\frac12\chi(\R A)+\chi_c(\R P_\e-\R A)$ and
$\chi_c(\R P_\e-\R A)=\chi_c(\Int \R P_\e)-\a_+= \chi(\Int \R
P_\e)-\a_+$, which yields finally $$ \align \D_\e&=b_2(\C
P)+(g_a(\C A)-1+ \frac12\chi(\R
A)-\frac12\mu+r)-\nu-n+n_\om+\chi(\Int \R P_\e)- \a_+\\ &=b_2(\C
P)-\nu+g(\C A)+\bra-(\a_{\roman{\Im}}+\a_+)-b_1(\Int \R P_\e)
+n_\om+b_2(\R P_\e)
\endalign
$$ \qed\enddemo

\subheading{4.4. Generalized Petrovskii inequalities} If we omit
in the Arnold-Viro-type inequalities any assumptions on the
singularities of an irreducible even curve $\C A$ in a
non-singular surface $\C P$, then some weaker estimates,
$b_2^\pm(\bar X^\circ)\ge0$, still can be used. One of them,
$b_2^+(\bar X^\circ)\ge0$, gives an estimate for $\frak p$ $$
b_2^+(\C P)+\frac12L.(K+L)+t_2+
%\min(\b,\til b_0(\C A'))
\frac12\min(\tilde b_0(\C A''),\tilde b_0(\C A)+2\b,b_2(\C A)-\nu)
\ge \frak p $$ (the author cannot say much about its use and
novelty). The other one, however, being the Petrovskii-type
inequality, contains a somewhat non-trivial information about
$\chi(\R X^\e)$ $$ \multline b_2^-(\C
P)+\frac12L.(K+3L)-\frac12\mu^-
%+\min(0,\til b_0(\C A')-\b)
+t_2-\b\\ +\frac12\min(\tilde b_0(\C A''),\tilde b_0(\C A)
+2\b,b_2(\C A)-\nu) \ge-\frac12\chi(\R X^\e),\\
\endmultline\tag4-31
$$ Applying (4-31) both for $\R X^+$ and $\R X^-$, and observing
that $\chi(\R X^\e)=\chi(\R P_\e)+\chi(\R P)-\chi(\R P_{-\e})$, we
obtain $$ \multline |\chi(\R P_+)-\chi(\R P_-)|\le2b_2^-(\C
P)+\chi(\R P)+L.(K+3L)-\mu^-+2t_2\\ +\min(\tilde b_0(\C
A'')-2\b,\tilde b_0(\C A),b_2(\C A)-\nu-2\b)
\endmultline
$$ In the classical case $\C P=\Cp2$, it gives $$ |\chi(\R
P_+)-\chi(\R P_-)|\le3k(k-1)+1-\mu^- +\min(0,\tilde b_0(\C
A'')-2\b,b_2(\C A)-\nu-2\b) \tag4-32$$ which is a refinement of
the generalized Petrovskii inequality stated by O. Ya. Viro
\cite{V1}.

\rk{Remarks} \roster
\item
Another version of the generalized Petrovskii inequalities for
real surfaces can be obtained if we use the estimate (4-19): $$
h^{1,1}(\C X^{\res})-2-2b_2(\bar E)+\hat\chi(\R E)\ge-\chi(\R X).
$$
\item
Recall that V.M. Kharlamov announced \cite{Kh1} a different sort
of the generalized Petrovskii inequalities. The relation of
Kharlamov's generalization to the generalization of Viro and to
the formula (4-32) seems to be an open question yet.
\endroster
\endrk

%%%%%%%%%%%%%%%%%%%%%%%%%%%%%%%%%%%%%%%%%%
\heading {\S5. Computation of the forms $\q_\pm$}
\endheading
\subheading{5.1. The local partition forms} Assume  that
$f\:(\C^2,0)\to(\C,0)$ is an isolated singularity defined by a
real polynomial $f(x,y)$ with the zero locus $\C A\subset\C^2$. By
\AG-morsification of $f$ (``\AG'' refers to A'Campo and
Gusein-Zade \cite{AC,GZ}) we mean a small real deformation,
$f^\tau$, of $f$ having $\mu$ non-degenerate {\it real} critical
points in an $\e$-ball, $\R B\subset\C^2$,
$0\!<\!\!<\!|\tau|\!<\!\!<\!\e\!<\!\!<\!1$, around zero ($\mu$ is
the Milnor number) and the maximal possible number,
$\frac12(\mu+\rho-1)$, of saddle points for $f^\tau|_{\R B}$,
which all lie on the same level curve $\C A^\tau=\{f^\tau=0\}$
($\rho$ here, like in \S4, is the number of the branches of $\R A$
at $0$).
Along with such objects as $\R B_\pm$, $\Cal H=\Cal H_+\oplus \Cal
H_-$, which were assigned to a singularity $f$ in 2.4, we consider
their deformations, $\R B_\pm^\tau=\R B\cap\{\pm
f^\tau(x,y)\ge0\}$, and ${\Cal H}^\tau_\pm=H^0(\R
B^{\tau}_\pm\setminus\R A^\tau)$, $\Cal H^{\tau}=\Cal
H_+^{\tau}\oplus \Cal H_-^{\tau}$. Denote by $W_i^\tau$ the
closures of the connected components of $\R B^{\tau}\setminus\R
A^\tau$, so that the negative indices, $i=-l,\dots,-1$ are used
for those of the components which lie in the interior  of $\R B$,
and positive, $i=1,\dots,2\rho$, for those which have common
points with $\R S=\d(\R B)$, if $\rho\ge1$. In the case of a dot
singularity (i.e., $\rho=0$), there is only one component bounded
by $\R S=\d\,\R B$, which we denote by $W^\tau_0$.

The regions $W^\tau_i$ associated to $f^\tau$ will be called {\it
\AG-regions} and the diagram $(\R B,\R B\cap\R A^\tau)$
characterizing the mutual adjacencies of $W^\tau_i$, will be
called {\it \AG-diagram}.

Denote by $W_i$ the closures of the connected components of $\R
B\setminus\R A$ which are deformed into $W_i^\tau$, $i\ge0$. Like
before, we put $\sign(W_i^\tau)\in\{+,-\}$ for the sign of
$f^\tau$ inside $W^\tau_i$, and $\sign(W_i)\in\{+,-\}$ for the
sign of $f$ inside $W_i$. We denote by $w_i^\tau\in\Cal H^\tau$
the generators representing $W_i^\tau$ and by $w_i\in H^0(\R
B\setminus\R A)\cong\Cal H$ the generators representing $W_i$.
Let us define a quadratic form, $\q^{\tau}\:{\Cal H}^{\tau}\to\Q$,
putting $$ \align
%  1
\q^{\tau}(w^\tau_i,w^\tau_j)&=0
\phantom{AAAAAAAAAAAAAAAAAAAAAAAA\,}\ \text{ if }
\sign(W_i^{\tau})\ne\sign(W_j^{\tau})
\\
%  2
 \q^{\tau}(w^\tau_i,w^\tau_j)&=\frac12\ord( W^\tau_i\cap  W^\tau_j);
\phantom{AAAAAAAAAAAAAAA}\text{ if }
\sign(W_i^{\tau})=\sign(W_j^{\tau}), \ i\ne j\\
%  3
 \q^{\tau}(w^\tau_i,w^\tau_i)&=
\frac12\ord( W^\tau_i\cap\Cl(\R B_\e^\tau-W^\tau_i))-
2\chi(W^\tau_i\setminus\R S) \text{\phantom{AAAA}\, if $-l\le i\le
2\rho$}\\
\endalign
$$
\comment
%  3
 \q^{\tau}(w^\tau_i,w^\tau_i)&=
\frac12\ord( W^\tau_i\cap\Cl(\R B_\e^\tau-W^\tau_i))
\phantom{AAAAAAA}\, \text{if \ } 1\le i\le2\rho,\ \text{for\ }
\rho\ge1\\
%  4
 \q^{\tau}(w^\tau_0,w^\tau_0)&=
\frac12\ord( W^\tau_0\cap\Cl(\R B_\e^\tau-W^\tau_0))
-2\chi_c(W^\tau_0) \text{\phantom{A} if $\rho=0$}\\
\endalign
$$
\endcomment
where $\ord$ denotes the number of points. Let $\Cal E$ denote the
subspace of ${\Cal H}^\tau$ generated by $w_i^\tau$,
$i\in\{-1,\dots,-l\}$. If the restriction of $\q^{\tau}$ to $\Cal
E$ is non-degenerated, then we obtain a direct sum decomposition,
$\Cal H^\tau=\Cal E\oplus \Cal E^\bot$, where $\Cal E^\bot$
denotes the orthogonal complement to $\Cal E$ with respect to the
quadratic form $\q^{\tau}$. Denote by $\hat w_i^\tau$, the
component of $w_i^\tau$ in $\Cal E^\bot$ with respect to this
direct sum decomposition, where $i=1,\dots,2\rho$ for $\rho\ge1$,
and $i=0$ for $\rho=0$.

Consider $\Cal E_\pm=\Cal E\cap\Cal H^{\tau}_\pm$ and $\Cal
E_\pm^\bot=\Cal E^\bot\cap\Cal H^{\tau}_\pm$. It is obvious that
$\Cal E^\bot=\Cal E_+^\bot\oplus \Cal E_-^\bot$ and that $\Cal
E_\pm^\bot$ is the orthogonal complement to $\Cal E_\pm$ with
respect to the form $\q^{\tau}_\pm=\q^{\tau}|_{{\Cal
H}^\tau_\pm}$, provided the latter is non-degenerated.

\tm{5.1.1. Theorem} Assume that $f\:(\C^2,0)\to(\C,0)$ is an
isolated real singularity and $f^\tau$ its \AG-morsification. Then
\roster
\item
$f$ is an \sQQ-singularity if and only if the restriction of
$\q^{\tau}$ to $\Cal E_+$ is non-degenerated.
\item
Assume that $f$ is an \sQQ-singularity. Then
$\bar{\q}(w_i,w_j)=2\q^{\tau}_+(\hat w^\tau_i,\hat w^\tau_j)$, for
any \AG-regions, $W_i$, $W_j$, if $0\le i,j\le2\rho$,
$\sign(W_i)=\sign(W_j)=+$.
\endroster
\endtm

Replacing $f$ by $-f$, we obtain a version of Theorem 5.1.1 for
$\q^{\tau}_-$ and $\Cal E_-$. Combined together, these two
versions imply

\tm{5.1.2. Theorem} Assume that $f^\tau$ is like in Theorem 5.1.1.
Then \roster
\item
$f$ is an \sQR-singularity if and only if the restriction of
$\q^{\tau}$ to $\Cal E$ is non-degenerated.
\item
Assume that $f$ is an \sQR-singularity. Then
$\q(w_i,w_j)=\q^{\tau}(\hat w^\tau_i,\hat w^\tau_j)$, for any
\AG-regions, $W_i,W_j$, $0\le i,j\le2\rho$.
\endroster
\endtm

Note the above theorems reduce the problem of calculating the
forms $\q$, $\q_\pm$ to an elementary combinatorial analysis of
\AG-diagrams (whose construction is known due to \cite{AC,GZ}) and
some trivial linear algebra (completion the squares of $w^\tau_i$,
$-l\le i\le -1$, in the form $\q^{\tau}$).

%%%%%%%%%%%%%%%%%%%%%%%%%%%%%%%%%%%%%%%%%%
\subheading{5.2. Contraction of $2$-dimensional polyhedra in
rational homology manifolds}
Assume that $Z$ is a compact oriented $\Q$-homology 4-manifold,
for simplicity, a polyhedron (for our purpose, it suffices to
consider only Whitney stratified pseudo-manifolds, which are known
to carry a polyhedral structure), and $K\subset Z-\d Z$ a
2-dimensional sub-polyhedron. Assume that $H_1(K;\Q)=0$, the
inclusion homomorphism, $\roman{in}\:H_2(K;\Q)\to H_2(Z;\Q)$, is
monomorphic and the restriction of the intersection form in $Z$ to
$\Cal E=\roman{in}(H_2(K;\Q))$ is non-degenerated. Analyzing the
long homology sequence of the pair $(Z,K)$, one easily obtains
that the quotient space $Z/K$ is also a $\Q$-homology manifold,
whose intersection form is isomorphic to the restriction of the
intersection form in $Z$ to the orthogonal complement, $\Cal
E^\bot$, of $\Cal E$ in $H_2(Z;\Q)$.

%%%%%%%%%%%
%%%%  5.3       5.3 5.3 5.3 5.3 5.3
%%%%%%%%%%%

\subheading{5.3. Proof of Theorems 5.1.1} Let us denote by $\C X$
and $\C X^\tau$ the affine surfaces defined in $\C^3$ by the
equations $f(x,y)-z^2=0$ and $f^\tau(x,y)-z^2=0$, and by $\pi$,
$\pi^\tau$ the projections of $\C X$, $\C X^\tau$ to the
$(x,y)$-plane, $\C^2$. Put $\C U=\pi^{-1}(\C B)$, $\C
U^\tau=(\pi^\tau)^{-1}(\C B)$, like in 2.4, and let
${\G}_i^\tau={\pi}^{-1}(W_i^\tau)$, $i\in\{-l,\dots,2\rho\}$.
Attaching 2-handles to $\overline U^{\tau}$ along the link
$L=\d\,\R U^{\tau}$, with the canonical framing, one obtains a
$4$-manifold, which we denote by $\overline U_L^{\tau}$. Similarly
attaching $2$-handles to $\overline U$, we obtain a $\Q$-homology
$4$-manifold, $\overline U_L$. In the case  $\rho\ge1$, we denote
by $\hat{\G}_i^\tau\subset \overline U_L^{\tau}$,
$i=1,\dots,2\rho$ the union of ${\G}_i^\tau$ with the core of the
corresponding handle. In the case of a positive dot singularity,
let $\hat{\G}_0^\tau$ be the union of ${\G}_0^\tau$ with the cores
of the both handles of $\overline U_L^{\tau}$.

The fundamental classes, $[{\G}_i^\tau]$, $i=-1,\dots,-l$, form a
basis in $H_2(\C U^\tau)$. Similarly, those of these classes, with
$\sign(W_i^\tau)=+$, form a basis of $H_2(\overline U^\tau)$ (the
both facts follow, for instance, from the description of the
Milnor form in \cite{AVG, section 1.4}). Together with
$[\hat{\G}_i^\tau]$, $i\ge0$, $\sign(W_i)=+$, the latter classes
form a basis of $H_2(\overline U^\tau_L)$. Since $\C U^\tau$ is
homotopy equivalent to a wedge of $2$-spheres, as was mentioned in
3.3, the set $R=\cup_{-l\le i\le-1}\G_i^\tau$ is a spine (a
regular deformational retract) of $\C U^\tau$, whereas
$R_+=R/\conj$ is a spine of $\overline U^{\tau}$. This implies
that $\C U\cong\C U^\tau/R$ and $\overline U\cong\overline
U^{\tau}/R_+$.

The homology sequence of the pair $(\overline
U^{\tau},\d\,\overline U^{\tau})$, implies that the intersection
form on $\overline U^{\tau}$ is non-degenerated if and only if
$\d\,\overline U^{\tau}$ is a $\Q$-homology sphere. On the other
hand, it is known that the intersection form in $\C U^{\tau}$ is
described by the restriction $\q^{\tau}|_{\Cal E}$ (see \cite{AVG,
1.4}). Together with Proposition 2.1.1, this implies that the
intersection form in $\overline U^{\tau}$ is described by the
restriction of $\bar q^{\tau}=2\q^{\tau}$ to ${\Cal E_+}$ and
proves the first part of the Theorem.

Similar calculations show that the formulae for the intersection
form in $H_2(\overline U_L^{\tau})\cong \Cal H_+^{\tau}$ are
almost the same as for $2\q^\tau_+$. The distinction arises only
for the self-intersections, $\la\hat{\G}_i^\tau,
\hat{\G}_i^\tau\ra_{\overline U_L^{\tau}}$, $i=1,\dots,2\rho$, for
$\rho\ge1$, which are equal to
$2q^\tau_+(w_i,w_i)-2\chi(W_i^\tau\cap\R S)$. On the other hand,
it follows from the definition of $\l$ and $\q$ that $$
\bar{\l}([\d {\G}_i],[\d {\G}_j])= \la\hat{\G}_i,
\hat{\G}_j\ra_{\overline U_L}= 2q_+(w_i,w_j)-2\chi(W_i\cap\R
S)\delta_{ij} $$ where $\delta_{ij}$ is the Kronecker symbol.
Since, obviously, $\chi(W_i\cap\R S)=\chi(W_i^\tau\cap\R S)=1$,
for $\rho\ge1$,
% (cf. \cite{AC, GZ}, or \cite{AVG, 1.4}),
and $\chi(W_i\cap\R S)=0$ for $\rho=0$, we obtain the second part
of Theorem 5.1.1 applying the observation in section 5.2 to
$Z=\overline U_L^{\tau}$ and $K=R_+$. \qed

\rk{Remarks} \roster\item Calculation of $\la\hat{\G}_i,
\hat{\G}_j\ra_{\overline U_L}$ is quite elementary: we just
combine the formula (4-1) with the Example in section 2.4.3
showing that the form $\overline \q_\pm$ for a cross-like node is
described by the matrix $\pmatrix\frac12&\frac12\\
\frac12&\frac12\endpmatrix$.
\item
An elementary analysis of the links at the points
$x\in\G_i\cap\G_j$ shows that $\overline U^\tau$ (and thus
$\overline U_L^\tau$) is a topological $4$-manifold. After it is
naturally smoothed, the intersection of $\G_i$ and $\G_j$ in
$\overline U^\tau$ becomes transversal (see \cite{F2}).
\endroster
\endrk

\subheading{5.4. Local complex intersection forms for simple
singularities} As an illustration, we present below the matrices
$M$ of the forms $\q_+$ for the simple real surface singularities,
which can be easily computed applying the algorithm described
above to the \AG-diagrams (sketched in \cite{AC}, \cite{GZ}) of
real simple singularities.

In the matrices $M$, characterizing the singularities
$D_{2n+2}^-$, $D_{2n+3}^-$ and $E_7$, the lesser diagonal entry
corresponds to the most narrow region (the one bounded by the real
branches which form the angle $0$). The shape of the other
singularities is symmetric, so we omit the correspondence  between
the regions and the entries of the matrices.

$$ \align
%%%%%%%%%%%%%%%%%%%%%%%%      A-(2n-1)---        A-(2n-1)----
%A_{2n-1}:\phantom{A}\,&
f(x,y)=&\cases -x^{2n}+y^2 &\phantom{AA}\,\,
A_{2n-1}^-\phantom{AA} M=\phantom{AAA\,} \pmatrix \frac n2&\frac
n2\\ \frac n2&\frac n2
\endpmatrix
\\
%%%%%%%%%%%%%%%%%%%%%%%     A-(2n-1)+++        A-(2n-1)+++
x^{2n}-y^2; &\phantom{AAA} A_{2n-1}^+\phantom{AA} M=\phantom{a}\,
\pmatrix \frac{2n-1}{2n}&\frac1{2n}\\ \frac1{2n}&\frac{2n-1}{2n}
\endpmatrix
\\
%%%%%%%%%%%%%%%%%%%%%%%%      A-(2n-1)000---        A-(2n-1)000---
x^{2n}+y^2; &\phantom{AA}\,\, A_{2n-1}^{\circ-}\phantom{AA}
M=\phantom{AAAA} \pmatrix2n-2\endpmatrix
\\
%%%%%%%%%%%%%%%%%%%%%%%%%%%%%%%%%%%%%%%%%%%%
%%%%%%%%%%%%%%%%%%%%%%%     A-(2n-1)000+++        A-(2n-1)000+++
-x^{2n}-y^2; &\phantom{AA}\,\, A_{2n-1}^{\circ+}\phantom{AA}
M=\phantom{AAAAA\,\,} \pmatrix0\endpmatrix
\endcases
\\
%%%%%%%%%%%%%%%%%%%%%%%%%%%%%%%%%%%%%%%%5
%\align
%A_{2n}:\phantom{AA}&
f(x,y)= &\cases \pm x^{2n+1}+y^2; &\phantom{A}\,\,
A_{2n}^-\phantom{AAA}\, M=\phantom{AAAAA\,} \pmatrix 2n
\endpmatrix
\\
%%%%%%%%%%%%%%%%%%%%%%%%%%%      A-(2n)+++      A-(2n)+++
\pm x^{2n+1}-y^2; &\phantom{A}\,\, A_{2n}^+\phantom{AAA}\,
M=\phantom{AAAA}\,\, \pmatrix\frac{2n}{2n+1}\endpmatrix
\\
\endcases
%%%%%%%%%%%%%%%%%%  DDDDDD      DDDDDDDDDD  DDDDDD
\\
%\endalign
%$$
%$$
%\align
%%%%%%%%%%%%%%%%%%%%%%%%%%%%%%%%%%  D-(2n+2)     D-(2n+2)
f(x,y)= &\cases \pm x(x^{2n}-y^2); &\phantom{A}
D_{2n+2}^-\phantom{AA} M=\phantom{A}\, \pmatrix 1 &\frac12
&\frac12\\ \frac12 &\frac{n+1}2 &\frac n2\\ \frac12 &\frac n2
&\frac{n+1}2
\endpmatrix
\\
\pm x(x^{2n}+y^2); &\phantom{A} D_{2n+2}^+\phantom{AA}
M=\phantom{AAAAA\,} \pmatrix 2n\endpmatrix
\\
\endcases
\\
%%%%%%%%%%%%%%%   D-(2n+3)  D-(2n+3)    D-(2n+3)
%D_{2n+3}:\ \ &
f(x,y)= &\cases x(x^{2n+1}\pm y^2); &\!D_{2n+3}^-\phantom{AA}
M=\phantom{AA}\, \pmatrix 2n+1 &1 \\ 1 &1
\endpmatrix
\\
-x(x^{2n+1}\pm y^2); &\!D_{2n+3}^+\phantom{AA} M=\phantom{AA}
\pmatrix \frac{2n+3}4&\frac{2n+1}4 \\ \frac{2n+1}4 & \frac{2n+3}4
\endpmatrix
\endcases
\\
%%%%%%%%%%%%%%%%%%%%%%%%%%%%%   E-6     E-6      E-6
%E_6:\phantom{AAA} &
f(x,y)= &\cases x^4\pm y^3; &\phantom{AAAA} E_6^-\phantom{AAA}\,\,
M=\phantom{AAAAAA}\pmatrix6\endpmatrix
\\
-x^4\pm y^3; &\phantom{AAAA} E_6^+\phantom{AAA}\,\,
M=\phantom{AAAAAA} \pmatrix2\endpmatrix
\\
%%%%%%%%%%%%%%%%%%%%%%%%%%%%%%%%%%%%%%%%%%%%%%%%%%%%%%%%%%%%%%%%%%%%%%%%%%
\endcases
\\
%%%%%%%%%%%%%%%%%%%%%%%%%%%%%%%%%%%%%%%%   E-7     E-7
%E_7:\phantom{AAA} &
f(x,y)=
%\ \
&\pm y(x^3\pm y^2); \phantom{AAAA}\, E_7\phantom{AAAA}
M=\phantom{AAAA} \pmatrix \frac72 &\frac32 \\ \frac32 &\frac32
\endpmatrix
\\
%%%%%%%%%%%%%%%%%%%%%%%%%%%%%%%%%%%%%%     E-8        E-8
%E_8:\phantom{AAA} &
f(x,y)=
%\ \
&\pm x^5\pm y^3; \phantom{AAAAAA} E_8\phantom{AAAA}
M=\phantom{AAAAAA}\pmatrix8\endpmatrix
\endalign
$$

\subheading{5.5. Some other methods and examples of computation of
the forms $q_{x}$} The forms $q_{x}$ for real isolated normal
surface singularities, can be described in terms of the resolution
of singularities, which is, in principle, more general and
sometimes more convenient then using \AG-morsifications. We sketch
a method how to do it, although under certain restrictions on a
singularity yet. Consider a compact $\conj$-invariant regular
cone-like neighborhood, $\C U$ of a real surface \QR-point, $x$,
and denote by $D_i$, $i=1,\dots,\rho$,  the closures of the
connected components of $\R U-\{x\}$; $D_i$ are the topological
discs bounded by the components $L_1,\dots,L_\rho$ of the real
link $L=\d (\R U)$. Assume that the exceptional divisor, $\C E$,
of a very good resolution, $r\:\C U^{\res}\to\C U$, contains at
least one real (that is $\conj$-invariant) component. Then
non-real components do not intersect $\R U^{\res}$, since $\C E$
is connected and does not contain triple intersection points.
Consider an intermediate resolution, $\res'\:\C U'\to\C U$, where
$\C U'$ is obtained from $\C U^{\res}$ by contraction of all the
non-real components of $\C E$. The following property is satisfied
in numerous examples: there exists a disjoint set, $\C
E^1,\dots,\C E^n\subset \C U'$, of the components of the
exceptional divisor of $\res'$ with the 2-sided real parts $\R
E^i\subset \R U'$, so that $\R U'-\bigcup_{i=1}^n \R E^i$ split
into $\rho$ orientable connected components, and for the closures
of these components, which we denote by $F_i$, $\res'(F_i)=D_i$,
$i=1,\dots,\rho$. Choose arbitrarily orientations of $D_i$ and the
coherent orientations of $F_i$. Put  $\e_k=1$, $k=1,\dots,n$, if
the surfaces, $F_i$, $F_j$, adjacent to $\R E^k$ from its two
sides (possibly $i=j$) induce the opposite orientations on $\R
E^k$, otherwise put $\e_k=-1$. Denote by $\a_k$ the
self-intersection of $\C E^k$ in $\C U'$. Note that $\a_k$ are
determined by the resolution graph of $\C U^{\res}$ and can be
calculated by the well-known continued fractions algorithm
\cite{HKK} (this algorithm, in turn, is justified by the remark in
5.2).

Then for $\gamma_i=[L_i]\in H_1(L)$ we have $$\align
\l_x(\gamma_i,\gamma_j)&=-\frac14 \sum_{\R E^k\subset F_i\cap
F_j}\e_{k}\a_k, \text{\ \ for $i\ne j$}\\
\l_x(\gamma_i,\gamma_i)&=\frac14 \sum_{\R
E^k\subset\d\,F_i}\e_{k}\a_k-\chi(F_i)\\
\endalign$$
The proof is essentially reduced to a quite elementary analysis of
the example 5.5.1 below.

Consider now the case of $\C E$ not containing $\conj$-invariant
components. Then $\R U^{\res}=\R U$ is a non-singular surface,
which is homeomorphic to a wedge of $\rho$ discs, hence $\rho=1$.
Furthermore, $\C E-\{x\}$ splits into a pair of connected
components. We denote by $E'$, $E''$ the closures of these
components and note that $E'$ and $E''$ are permuted by $\conj$
and homeomorphic to $\overline E$. The definition of the
resolution graph admits an obvious extension to the quotients,
$\overline E\subset\overline U^{\res}$; namely, such a graph
(being a tree in our case) characterizes likewise the
intersections of the components of $\overline E$ in $\overline
U^{\res}$ (one can call it a {\it quotient resolution graph}).
Such a graph almost coincides with a ``half'' of the usual
resolution graph, i.e., a subgraph representing the components of
$E'$ (or, equally, of $E''$). The distinction arises only with the
weight of one vertex corresponding to the component, $E_x$, of
$E'$, which contains $x$ (the weight in the quotient resolution
graph is less by $1$).

Let $\hat U$ be obtained from $\overline U^{\res}$ by contraction
of all the components of $\overline E$ except $E_x$; then $\hat U$
is a $\Q$-homology manifold containing $E_x$ as a deformational
retract. Using the continued fractions algorithm, we can determine
the self-intersection $\a=\la E_x,E_x\ra_{\hat U}$ and obtain,
applying the remark in 5.2, that $$\l_x([L],[L])=-\frac4{\a}-1$$
where $4$ appears as the square of $\la E_x,\R U^{\res}\ra_{\hat
U}=\pm2$.

{\bf Example 5.5.1.} Consider a rational $(-m)$-curve $\C E$,
$m>0$, $m\in\Q$,
 in a \QR-surface $\C X'$ such that
the real part $\R E$ of $\C E$ (but not necessarily the whole $\C
E$) is smooth. Then the linking form $\l_{x}$ of the singularity
which appears in $\C X=\C X'/\C E$ after contraction of $\C E$, is
described by the matrix $(-m)$, if $\R E$ is one-sided  in $\R X$,
and by the matrix $\pmatrix -\frac{m}4&\frac m4\\ \frac m4&-\frac
m4
\endpmatrix
$
if it is two-sided.

We can apply this method to calculate, for instance, the form
$q_{x}$ for the singularity $A_{2n-1}^+$, since the latter appears
after
 contraction of a rational curve $\C E$  with a pair of complex-conjugated
imaginary singularities of the type $A_{n-1}$ and with $\la \C E,
\C E\ra_{\C X'}=-\frac2n$. This gives a matrix $\pmatrix
-\frac{1}{2n}&\frac1{2n}\\ \frac1{2n}&-\frac{1}{2n}
\endpmatrix$
of the form $\l_{x}$.

{\bf Example 5.5.2.} We call a real surface singularity {\it
quasi-cuspidal} if its real link $\R M$ has one component (an
obvious example is the suspension over an unibranch curve
singularity). If, for a quasi-cuspidal singularity, $\R U^{\res}$
is orientable, then, putting $\C U'=\C U^{\res}$, $F=\R U^{\res}$,
and applying the algorithm described in 5.4, we obtain
$\l([L],[L])=-\chi(F)=2g-1$, where $g$ is the genus of $F$. If we
are given an unibranch curve singularity, so that the suspension
surface singularity has orientable real part $\R U^{\res}$, then
the residue form is determined by the value $\q_+(v)=2g$, where
$v\in\Cal H_+\cong\Z$ is a generator.

Note that $\R U^{\res}$ is orientable if the real components of
$\C E$ have even self-intersections. The Lefschetz fixed point
formula for the involution $\conj$ implies moreover that $\chi(\R
U^{\res})=1-m$, where $m=2g$ is the number of real components of
$\C E$. For instance, for quasi-cuspidal singularities $A_{2n}^-$
and $D_{2n+2}^+$, the number of real components is equal to $2n$,
for $E_6^+$ and $E_8$, this number is $2$ and $8$ respectively,
which determine their forms $q_+$ (cf. the table in 5.4).

Another example is the singularity $f\:(\C^2,0)\to(\C,0)$,
$f(x,y)=x^{2n}\pm y^{2n-1}$, where $n$ is even. In this case $\C
U^{\res}$ is spin, and thus $\R U^{\res}$ is orientable.
Furthermore, it is not difficult to determine the number,
$m=4n-2$, of the real components of $\C E$, which implies that
 $q_+(v,v)=(4n-2)$.
It contrasts to the case of odd $n$, in which $\C U^{\res}$ is not
Spin and $q_+=0$ (see the following example).

{\bf Example 5.5.3.} For the singularity at $0\in\C^2$ defined by
$f(x,y)=x^{2n}\pm y^{2n-1}$, where $n$ is odd, the form $\q_+$
vanishes. To see it, we consider the double plane, $\C
X^-\to\Cp2$, branched along the projective closure, $\C A$, of the
curve $\{f=0\}$ defined by the equation $x^{2n}\pm
y^{2n-1}z-w^2=0$. The quotient ${\overline X}^-$ is a double
covering over $\Cp2/\conj\cong S^4$ branched along
$\A_+=\Rp2_+\cup\overline A$ (cf. \cite{Ar}), where $\overline A$
is a $2$-disc, because $\C A$ is a rational curve. Since $\Rp2_+$
is also a $2$-disc, we have $\A_+\cong S^2$ and ${\overline X}^-$
is a homotopy sphere, because it is simply connected, as $\C X$ is
(actually, one can show that ${\overline X}^-$ is diffeomorphic to
$S^4$). This implies that $\la \R X, \R X\ra_{\overline X}=0$ and
thus $q_+=0$, because $\R X$ is a torus with a unique singularity.

{\bf Example 5.5.4.} Assume that $p,q\ge1$, $p+q=2n$, $(p,2n)=1$
and $n$ is odd. Then the singularities defined at $0$ by
$f(x,y)=x^{2n}-y^{p}$ and $g(x,y)=x^{2n}-y^{q}$, have the opposite
forms, i.e., $q_{f}(v_f,v_f)=-q_{g}(v_g,v_g)$, where $v_f\in\Cal
H_{f+}\cong\Z$, $v_g\in\Cal H_{g+}\cong\Z$ are generators.

To prove it, we apply the same arguments as in the previous
example to the curve $\C A=\{x^{2n}-y^{p}z^{2n-p}=0\}\subset\Cp2$.

\heading {\S6. Some examples and applications}
\endheading
\subheading{6.1. Arrangements of hyperplanes} Consider a
hypersurface $\C A\subset\Cp{d}$ which splits into $m=2k$ real
hyperplanes in generic position; then $\C A$ is a \sQ-hypersurface
(see Example 2.4.3 in section 2.4). Let $W_1,\dots,W_m$ denote the
partition regions (i.e., the polyhedra bounded by the hyperplanes)
of $\C A$. If $d-k$ is even, then we can choose as $\Omega$ (cf.
section 2.5) one of the hyperplanes; if $d-k$ is odd, then we put
$\Omega=\oo$.

We assume first that $\G=W_i\cap W_j$ is a {\it connected}
polyhedron (which is always the case if $m$ is sufficiently
large). Let $f_{\Gamma}(t)$ be the face-counting polynomial
defined as $$ f_{\Gamma}(t)=\sum_{0\le k\le s}f_k t^k $$ where
$s=\dim\G$ and $f_k$ is the number of $k$-faces of $\Gamma$.
If $\Gamma\subset\Omega$ and $\Omega$ locally (near $\G$)
separates $W_i$ and $W_j$, then we put $\e(W_i,W_j|\Omega)=-1$,
otherwise, we put $\e(W_i,W_j|\Omega)=1$. Split the regions $W_i$
into positive and negative assigning a sign,
$\sign(W_i)\in\{+,-\}$, so that $W_i$ and $W_j$ have no common
$(d-1)$ face if $\sign(W_i)=\sign(W_j)$. Put
$\sign(W_i,W_j)=\frac12(\sign(W_i)-\sign(W_j))$.

The formulae (2-9) easily imply the following result

\tm{6.1.1. Theorem} Assume that $\G=W_i\cap W_j$ is a connected
polyhedron of dimension $s$. Then $$ \align
\phi(w_i,w_j|\Omega)&=\e(W_i,W_j|\Omega)
(-1)^{\frac12(d^2-s+\sign(W_1,W_2))}2^{1-d}f_{\G}(-2)\\
\endalign
$$ for any $i,j\in\{1,\dots,m\}$. \qed\endtm

%In particular, if $W_i$, $W_j$ are the both positive or the both negative,
%then $\phi(w_i,w_j|\Omega)=\e(W_i,W_j|\Omega)
%(-1)^{\frac12(d^2-s)}2^{1-d}f_{\G}(-2)$,
 Recall that one can associate certain real quasi-smooth toric variety,
$\C T_\G$, to a polyhedron $\G$. The algebraic structure of $\C
T_\G$ depends on the geometry of $\G$ (which actually must have
rational vertices). However, the Poincare polynomial, $\roman
P_{\C T_\G}(t)$ of $\C T_\G$, is determined by the combinatorial
type of $\G$, namely, $\roman P_{\C T_\G}(t)=f_{\G}(t^2-1)$ (see
\cite{Da}). In particular, this implies that $$f_\Gamma(-2)=\roman
P_{\C T_\G}(i)=\s(\C T_\G)=-\chi(\R T_\G)$$
{\bf Example:} If $\G=W_i\cap W_j$ is a simplex of dimension $s$,
then $\roman P_{\C T_\G}=\roman P_{\Cp{s}}$,
 so $|\phi(W_i,W_j|\Omega)|=(\frac12)^{d-1}$
if $s$ is even and $\phi(W_i,W_j|\Omega)= 0$ if odd. In the case
of a prism, $W_i\cap W_j\cong\D\times[0,1]$, with any convex
polyhedron $\D$ as the base, we have $\phi(W_i,W_j)=0$. More
generally, one can use the obvious relation
$f_{\G_1\times\G_2}=f_{\G_1}\times f_{\G_2}$ for evaluation of
$\phi(W_i,W_j)$ if $W_i\cap W_j\cong \G_1\times\G_2$.

If $W_i\cap W_j=\G_1\cup\dots\cup\G_r$ is a union of several
connected polyhedra, then $\phi(w_i,w_j|\Omega)$ is the sum of the
expressions given by Theorem 6.1.1, calculated for each polyhedron
$\G_i$.

\tm{6.1.2. Theorem} Assume that $\C A\subset\Cp{d}$ is as above,
an arrangement of hyperplanes, $d=2n$ and $\pi\:X\to \Cp{d}$ the
double covering branched along $\C A$. Denote by $\C X^\e$, $\e=+$
or $\e=-$, the complex variety $X$ endowed with one of the two
real structures (defined by the complex conjugations) lifted from
$\Cp{d}$, and put $\R \roman P_\e=\pi(\R X^\e)$. Let
$\phi_{\ee}\:H^0(\R \roman P_\e-\R A)\to \Q$ denote the components
of the partition form, $\phi=\phi_+\oplus\phi_-$. Then $$ \align
\s_{-\ee}(\phi_{\ee})&=\frac12(b_d^{-\ee}(X)-\kapnew)\\
\s_{\ee}(\phi_{\ee})&=\frac12(b_d^{\ee}(X)+\kapnew+\chi(\R
X^\e))-1\\ \s_0(\phi_{\ee})&=\sum_{k=0}^{d-1}\binom{m-2}{k}
\endalign
$$ where $\ee$ and $\kapnew$ are like in subsection 3.1.
\endtm

This result shows that the generalized Arnold-Viro inequalities
are indeed equalities for arrangements of hyperplanes.

\demo{Proof} The Arnold-Viro inequalities (3-1)--(3-2) give $$
\aligned \s_{-\ee}(\phi_{\ee})+\s_0(\phi_{\ee})&\le
\frac12(b_d^{-\ee}(X)-\kapnew)+\delta(\C X^{-\e})\\
\s_{\ee}(\phi_{\ee})&\le \frac12(b_d^{\ee}(X)+\kapnew+\chi(\R
X^\e))-1\\
\endaligned\tag6-1
$$ and thus $\s_+(\phi_{\ee})+\s_-(\phi_{\ee})+\s_0(\phi_{\ee})
\le\frac12(b_d(X)+\chi(\R X^\e))-1+\delta(\C X^\e)$

Combining together these inequalities for $\e=+$ and $\e=-$, we
obtain $$ N^m_d\le b_d(X)+\frac12(\chi(\R X^+)+\chi(\R X^-))-2+
\delta(\C X^+)+\delta(\C X^-) \tag6-2 $$ where
$N^m_d=N^m_d=\sum_{k=0}^{d}\binom{m-1}{k}=\binom{m-2}{d}
+2\sum_{k=0}^{d-1}\binom{m-2}{k}$ is the number of the partition
regions for a generic arrangement of $m$ hyperplanes in $\Rp{d}$.
It is easy to check that $b_d(X)=\binom{m-2}{d}+1$ for even $d$
(while $b_d(X)=\binom{m-2}{d}$ for odd $d$). Now the obvious
relation $\chi(\R X^+)+\chi(\R X^-)=2\chi(\Rp{d})=2$ and the
estimate of $\delta(\C X^\pm)$ in the following proposition show
that (6-2) is indeed an equality and thus all the intermediate
estimates, including (6-1),
 are also equalities.

\tm{6.1.3. Proposition} $\delta(\C X^\pm)\le
N_{d-1}^{m-1}=\sum_{0\le k\le d-1}\binom{m-2}k$.
\endtm

The proof is given in subsection 7.4. \qed
\enddemo

\subheading{6.2. Arnold inequalities for quasi-cuspidal surfaces}
A real surface $\C X$ will be called {\it quasi-cuspidal} if the
links, $\R M_x$, for $x\in \R X$ are circles, i.e., $\R X$ is
topologically non-singular. We call a connected component
$\G\subset \R X$ {\it elliptic, parabolic or hyperbolic} if
$-\la\G,\G\ra_{\C X}$ is positive, zero or negative respectively.
Such terms are motivated by the relation $-\la\G,\G\ra_{\C
X}=\chi(\G)$ satisfied provided $\G\cap\Sing(\R X)=\oo$. More
generally, if $\G\cap\Sing(\R X)$ contains only \QR-singularities,
Theorem 2.2.4 shows that one can treat $-\la\G,\G\ra_{\C X}$ as
the {\it weighted Euler characteristic} of $\G$, in which points
$x\in\G\cap\Sing(\R X)$ are counted with the weights $-\l_{x}([\R
M_x],[\R M_x])$.

Given a real quasi-cuspidal \QQ-surface $\C X$, denote by $c^+$,
$c^0$, $c^-$ the number of {\it oriented} elliptic, parabolic and
hyperbolic components of $\R X$. Inequalities of Theorem 4.2.1
give estimates for $c^+$, $c^0$, $c^-$. (4-17) yields $c^-\le
p_g(\C X^{\res})$, so, for instance, $c^-=0$ for a singular
rational or Enriques surface. If $\C X$ is a quasi-cuspidal
surface with only $\Z/2$-inessential singularities, then Theorem
4.2.1 implies that $$ c^0+c^-\le\chi_a(\C X^{\res})+b_1(\C
X^{\res};\Z/2) $$ which gives for instance $c^0\le1$ for rational
surfaces and $c^0\le2$ for Enriques surfaces.
The estimates (4-19) and (4-20) can be used similarly and give an
information about $\chi(\R X)$, like in the non-singular case.

Note that none of the above estimate follows automatic from the
well known similar estimates for non-singular surfaces, because
$c^-$, $c^0$ (as well as $c^+$) may decrease (as well as increase)
after we pass to a resolution $\C X^{\res}\to\C X$. As a simplest
example, take a double plane, $\C X^+$, branched along a quartic
$\C A\subset\Cp2$, with $\R A$ being a single oval with $3$
ordinary cusps. $\R X^+$ has a parabolic component, which becomes
elliptic after resolution.

\rk{Remark} The generalized Arnold-Viro inequalities for
quasi-cuspidal surfaces have extremal properties analogous to the
well-known such properties for non-singular surfaces (see
\cite{Wi} \cite{V4}). Namely, the estimates (4-18), (4-20) can be
improved by $1$, unless all the components of $\R X$ are
orientable and parabolic. More generally, if (4-18) or (4-20) is
an equality for any \QQ-surface $\C X$, then there exists an
integral lifting, $a\in H_2(\R X)$, of the fundamental class $[\R
X]_2\in H_2(\R X;\Z/2)$, of $\R X$, such that the image of $a$ in
$H_2(\C X')$ (and thus, in $H_2(\C X)$) vanishes. This gives
certain restriction on the intersections (compare with the
analogous formulation in \cite{Kh2}, where the case of nodal
curves is considered).
\endrk

%%%%%%%%%%%%%%%%%%%%%%%%%%%%%%%%%%%%%%%%%%%%%%%%%%%%%
\subheading{6.3. Arnold inequalities for cuspidal and
quasi-cuspidal curves} A topologically non-singular complex curve
will be called {\it a cuspidal curve}; by a {\it quasi-cuspidal
curve} we mean a real curve, $\C A$, with the topologically
non-singular real part $\R A$. Accordingly, the singularities of a
cuspidal  curve will be called (generalized) {\it cusps}, and the
real singularities of a quasi-cuspidal curve, respectively, {\it
quasi-cusps}. The Smith inequality (7-2) implies that cusps are
\sQ-singularities (and, moreover, $\Z/2$-inessential
singularities).

Like for non-singular curves, the connected components of $\R A$
for a quasi-cuspidal real curve $\C A\subset\Cp2$ are all
null-homotopic in $\Rp2$ if $\C A$ has even degree $d=2k$.
Following the tradition, we call null-homotopic components {\it
ovals}. An oval is called even (or odd), if it lies inside an even
(respectively, odd) number of the other ovals.

Consider the double plane, $\pi\:X\to\Cp2$, branched along $\C A$
and let $\C X^\pm$ denote $X$ endowed with the complex conjugation
$\conj^\pm$, covering the complex conjugation in $\Cp2$, and put
$\Rp2_\pm=\pi(\R X^\pm)$. Like in the non-singular case, there
exists a single non-orientable partition component of $\Rp2-\R A$.
We denote it by $W_{\infty}$, and let
$\G_\infty=\pi^{-1}(W_\infty)$, assuming for definiteness that
$\sign(W_{\infty})=-$.
Given a \QQ-curve $\C A$, we denote by $p^+$, $p^-$ and $p^0$,
($n^+$, $n^-$ and $n^0$) the number of even (respectively, odd)
ovals, $C_i\subset\R A$, with the elliptic, parabolic and
hyperbolic component $\G_i\subset\R X^+$ (respectively,
$\G_i\subset\R X^-$), where $\G_i=\pi^{-1}(W_i)$ and
 $W_i$ is the partition region bounded by $C_i$ from outside.
We introduce furthermore the indicators,
$\e^+,\e^0,\e^-\in\{0,1\}$, showing the type of the components
$\G_{\infty}$. If $k$ is even and thus $\G_{\infty}$ is
non-orientable, we put $\e^++\e^0+\e^-=0$. If $k$ is odd, then
$\e^++\e^0+\e^-=1$ with the non-vanishing indicator being $\e^+$,
$\e^0$ or $\e^-$, if $\G_\infty$ is elliptic, parabolic, or
hyperbolic respectively. The terms $r$, $\nu$, $\mu^\pm$, are
defined as in \S4.

\tm{6.3.1. Theorem} Assume that $\C A\subset\Cp2$ is a real
quasi-cuspidal \QR-curve of degree $2k$. Then $$ \align
n^-+n^0+\e^-&\le
\frac12(k-1)(k-2)-{\frac12\mu^+}+\min(n^0,r-\nu-1)\\ p^-+p^0&\le
\frac12(k-1)(k-2)-{\frac12\mu^+}+\min(p^0,r-\nu
-\frac12(1-(-1)^k))\\ n-p^-&\le
\frac32k(k-1)-\frac12\mu^-+\min(p^0,r-\nu-\frac12(1-(-1)^k))\\
p-n^-+\e_+&\le \frac32k(k-1)-\frac12\mu^-+\min(n^0+1,r-\nu)
\endalign
$$
\endtm

In the case of cuspidal curves, we have $r=1$, $\nu=0$ and the
above inequalities are simplified as follows

\tm{6.3.2. Corollary} If $\C A\subset\Cp2$ is a real cuspidal
curve, then $$ \align n^-+n^0+\e^-&\le
\frac12(k-1)(k-2)-\frac12\mu^+\\
%+\min(0,n^0-\e_-)\\
p^-+p^0&\le
\frac12(k-1)(k-2)-\frac12\mu^++\frac12(1+(-1)^k)\tag6-3\\
n-p^-&\le \frac32k(k-1)-\frac12\mu^-+\frac12(1+(-1)^k)\tag6-4\\
p-n^-+\e^+&\le \frac32k(k-1)-\frac12\mu^-+1\\
%+\min(0,n^0-\e_+)\\
\endalign
$$
\endtm

\rk{Remarks} (1) The extremal properties of the Arnold
inequalities mentioned in 6.2 allow us, like in the case of
non-singular $\C A$, improve by $1$ some of the above estimates,
unless the curve $\R A$ has a rather special topology, with all
the components $\G\subset\R X^\pm$ being parabolic; for instance,
one can put $1$ instead of $\frac12(1-(-1)^k)$ in Theorem 6.3.1
and to drop $\frac12(1+(-1)^k)$ from (6-3) and (6-4).

(2) The term $\frac12(1+(-1)^k)$ can be omitted in (6-4) if $\C A$
is non-singular (cf. \cite{Wi,V4}), since in the extremal case we
get $p=p^0=n$ and (6-3) turns out to be stronger then (6-4).
Similar arguments show that we can also omit this term if
$\mu^--\mu^+\le2k^2-2$.
\endrk

\subheading{6.4. Pentics with three $A_4$-cusps} The following
example can give an idea how one can get an information about the
geometry of the singularities on a curve.
Assume that $\C A$ is a real pentic having three cusps, $x_i\in\R
A$, $i=1,2,3$, of the type $A_4$ (so that $\C A$ is rational).
Consider a line $\C L$ passing through the points $x_1$, $x_2$,
and let $y$ denote the third intersection point of $\R A$ with $\R
L$.

\tm{6.4.1. Proposition} The curve $\R A$ approaches the line $\R
L-\{y\}$ at points $x_1$, $x_2$ from the opposite sides (which
makes sense, because $\R L-\{y\}$ is bilateral in $\Rp2-\{y\}$).
\endtm

\demo{Proof} The partition form of a sextic $\C A'=\C A\cup\C L$
is singular by Corollary 4.3.3, since it is reducible. However, an
elementary combinatorial analysis of the possible mutual position
of the cusps and $\R L$ shows that this form cannot be singular if
$\R A$ comes to $\R L-\{y\}$ from the same side at $x_1$ and
$x_2$.
\enddemo

\rk{Remark} Proposition 6.4.1 can be easily derived also from
Gudkov's theorem characterizing possible mutual position of a
non-singular pentic and a line, however, our arguments are
essentially more elementary then the Gudkov's result.
\endrk

\heading \S7 Appendix A. The Smith theory estimates
\endheading
\subheading{7.1. The Smith sequence} Given  an involution on a
finite $CW$-complex, $c\:Z\to Z$, with the fixed point set $F$ and
the quotient $Z/c$, we can write a long homology {\it Smith exact
sequence} (with $\Z/2$-coefficients) $$ \dots\to H_{k+1}(Z/c,F)\to
H_k(Z/c,F)\oplus H_k(F)\to H_k(Z)\to H_k(Z/c,F)\to\dots $$ where
$H_k(F)\to H_k(Z)$ is the inclusion homomorphism, $H_k(Z)\to
H_k(Z/c,F)$ is induced by the quotient map $q\: Z\to Z/c$, and
 $H_{k+1}(Z/c,F)\to H_k(F)$ is taken from the homology sequence of the pair
$(Z/c,F)$ (see, e.g., \cite{Br, p.123}, or \cite{Wi, Appendix}).
Denote by $\nu_k$ the rank of the inclusion homomorphism
$H_k(F;\Z/2)\to H_k(Z/c;\Z/2)$.

Analyzing this sequence one easily obtains the following
estimates, called {\it the Smith inequalities}.

\tm{7.1.1. Theorem} For any $k\ge0$
$$b_k(Z/c,F;\Z/2)+b_k(F;\Z/2)\le
b_{k+1}(Z/c,F;\Z/2)+b_k(Z;\Z/2)\tag7-1$$ For any $l\ge0$ $$
b_l(Z/c,F;\Z/2)+\sum_{k\ge l}b_k(F;\Z/2)\le \sum_{k\ge
l}b_k(Z;\Z/2)\tag7-2 $$ and in particular $b_*(F;\Z/2)\le
b_*(Z;\Z/2)$, where $b_*$ stands for the sum of all Betti numbers.
$$ b_k(Z;\Z/2)\le b_{k}(Z/c;\Z/2)+b_k(Z/c,F;\Z/2)\tag7-3 $$ for
any $k\ge0$. The latter can be also formulated as $$
b_k(Z;\Z/2)\le
2b_k(Z/c;\Z/2)+b_{k-1}(F;\Z/2)-\nu_k-\nu_{k-1}\tag7-4 $$ where by
definition $\nu_{-1}=b_{-1}(\dots)=0$.
\endtm

Here (7-1) is proved obviously, adding (7-1) for all $k\ge l$, we
obtain
 (7-2).
(7-3) follows from that the dimension of the kernel of
$H_k(Z/c,F;\Z/2)\to {H}_{k-1}(F;\Z/2)$ and of the cokernel of
$H_{k+1}(Z/c,F;\Z/2)\to H_k(F;\Z/2)$ in the Smith sequence give in
sum $b_k(Z/c;\Z/2)$. (7-4) is obtained from (7-3) if we use the
expression of $b_k(Z/c,F;\Z/2)$ given by the exact sequence of the
pair $(Z/c,F)$.

\rk{Remarks} \roster
\item
Using augmented homology groups, we obtain the same inequalities
for $\til b_k$.
\item
Assume that $\dim Z=m$, $b_m(Z;\Z/2)=1$ and $\dim F\le m-2$. Then,
analyzing the first non-trivial terms in the Smith sequence, we
obtain $$ \align b_{m}(Z/c;\Z/2)&=1\\ b_{m-1}(Z/c,F;\Z/2)&\le
b_{m-1}(Z;\Z/2)+1\\ b_{m-1}(Z;\Z/2)&\le
b_{m-1}(Z/c;\Z/2)+b_{m-1}(Z/c,F;\Z/2)-1\\
&=2b_{m-1}(Z/c;\Z/2)+\til b_{m-2}(F;\Z/2)-1 \tag7-5\\
\endalign
$$ where (7-5) is an improvement of (7-3)--(7-4) in this special
case
\item
If $A\subset Z$ is a $c$-invariant subcomplex, then the Smith
sequence for the induced involution $Z/A\to Z/A$ gives $$ \dots\to
H_k(Z/c,F\cup A/c)\oplus H_k(F,F\cap A)\to H_k(Z, A)\to
H_k(Z/c,F\cup A/c)\to\dots $$ (all the groups are with
$\Z/2$-coefficients), which gives obvious relative versions of the
estimates (7-1)---(7-4).
\endroster
\endrk

\subheading{7.2. Proof of Lemma 3.1.3} The homology sequences of
$(\overline X,\R X)$ and the universal coefficients formula imply
$$ \delta(\C X)=b_{d+1}(\overline X,\R X)-b_{d+1}(\overline X)\le
b_{d+1}(\overline X,\R X;\Z/2)-b_{d+1}(\overline X) \tag7-6$$
Applying (7-2), with $l=d+1$, to the complex conjugation in $Z=\C
X$ and using (7-6) we obtain estimate (3-5). The next estimate,
(3-6), can be obtained from (3-5), since for a $d$-dimensional
complete intersection, $\C X$, with isolated singularities, we
have $b_k(\C X)=b_k(\Cp{d})$, for all $k$, except possibly $d$ and
$d+1$. The latter follows from the Lefschetz hyperplane section
theorem, homotopy equivalence of the Milnor fibers to wedges of
spheres, and from the exact sequence of the pair $(\C
X^\tau,U^\tau)$.

We may prove (3-7) in a slightly more general setting, with
$\Cp{d}$ being replaced by a real non-singular variety $\C P$ of
dimension $d$. Applying (7-3) to the deck transformation of the
double covering $p\: \C X\to \C P$ branched along a real reduced
hypersurface $\C A\subset \C P$, we obtain $$b_{k}(\C X;\Z/2)\le
b_{k}(\C P;\Z/2)+b_{k}(\C P,\C A;\Z/2) \tag7-7$$ which can be
improved by $1$ for $k=2d$ and $k=2d-1$ as (7-5) shows. Using the
improved inequality to estimate $b_{k}(\C X;\Z/2)$ in (3-5), we
obtain (3-7). Similarly, (7-4) combined with (3-6) gives (3-8).

%%%%%%%%%%%%%%%%%%%     7.3     7.3     7.3

\subheading{7.3. Proof of Lemma 4.1.4} The relations
$b_i(\overline X)=\frac12b_i(\C X)$, $i=1,3$, in (4-10) and (4-11)
follow, for instance, from that the mixed Hodge structure in
$H^i(\C X;\C)$ is pure. There is a more elementary way to derive
this relations: just to notice that $b_1(\C X)=b_1(\C X^\tau)$,
$b_1(\overline X)=b_1(\overline X^\tau)$, on one hand, whereas
$b_3(\C X)=b_3(\C X^{\res})$, $b_3(\overline X)=b_3(\overline
X^{\res})$, on the other hand, and then to use the relation
$b_i(\overline Z)=\frac12b_i(\C Z)$, which obviously holds for
non-singular real varieties, $\C Z$, and odd $i$ (cf. 8.1). The
relation (7-4) applied to the deck transformation of $\C X\to \C
P$ and its restriction to $\C X'$ and $\C X''$ gives $$ \align
b_1(\C X;\Z/2)&\le 2b_1(\C P;\Z/2)+\til b_0(\C A)\tag7-8\\ b_3(\C
X;\Z/2)&\le 2b_3(\C P;\Z/2)+b_2(\C A)-\nu\tag7-9\\ b_1(\C
X';\Z/2)&\le 2b_1(\C P;\Z/2)+\til b_0(\C A')\tag7-10\\ b_1(\C
X'';\Z/2)&\le 2b_1(\C P;\Z/2)+\til b_0(\C A'')\\
\endalign
$$ where
%$r=b_2(\C A;\Z/2)$ and
$\nu=\rank(H_2(\C A;\Z/2)\to H_2(\C P;\Z/2))$, like in 4.1. (7-8)
implies (4-10), since $b_1(\C X)\le b_1(\C X;\Z/2)$. (4-11)
contains $3+1$ estimates of $b_3(\overline X)$, and the first
among them follows from (7-9). The next two estimates follows from
$$ \align b_3(\overline X)&=b_3(\overline X',\d\overline
X')=b_1(\overline X') \le b_1(\C X')\le 2b_1(\C P;\Z/2)+\til
b_0(\C A')\\ 2b_3(\overline X)&=b_3(\C X)=b_3(\C X'',\d\C
X'')=b_1(\C X'') \le 2b_1(\C P;\Z/2)+\til b_0(\C A'')\\
\endalign
$$
%combined with (7-10).
The fourth estimate is a combination of (7-8) with the inequality
$b_3(\overline X)\le b_1(\overline X)+\b$, which follows from
(4-21).
The arguments in 7.2 give in the case $d=2$ $$ \align \delta(\C
X)&\le 2b_3(\C P;\Z/2)+b_2(\C A)-\nu-b_3(\overline X)\\ \delta(\C
X')&\le 2b_3(\C P';\Z/2)+b_2(\C A')-\nu'-b_3(\overline X')\\
\endalign
$$ Furthermore, for $\a_{\Im}\ge1$ we have $$ \align b_3(\overline
X')&=b_1(\overline X',\d\overline X')= b_1(\overline X)+\til
b_0(\d\overline X)=b_1(\overline X)+(\a_{\Im}-1)\\ b_3(\C
P';\Z/2)&=b_1(\C P',\d\,\C P';\Z/2)= b_1(\C P;\Z/2)+\til
b_0(\d\,\C P')= b_1(\C P;\Z/2)+(2\a_{\Im}-1)\\
\endalign
$$ which proves (4-12). (4-13) follows from (4-12) combined with
(4-7).

%%%%%%%%%%%%%%%%%%%%%%%%%%%%%%%%%%%%%
%%%%%%%%%%%%%%%%%%%%%%%%%%%%%%%%%

\subheading{7.4. Proof of Lemma 4.2.3} Let $\C U_{\Im}^{(2)}$
denote the union of the connected components of $\C U$ around the
$\Z/2$-essential imaginary singular points of $\C X$. Put $\C
X_{\Im}^{(2)}=\Cl(\C X-\C U_{\Im}^{(2)})$. We have $H_2(\overline
X';\Q)\cong H_2(\overline X_{\Im}^{(2)};\Q)$, since $\overline
X_{\Im}^{(2)}$ is obtained from $\overline X'$ by removing several
cones over $\Q$-homology spheres, and thus $\delta(\C
X')=\dim\ker(H_2(\R X;\Q)\to H_2(\overline X_{\Im}^{(2)};\Q))=
b_3(\overline X_{\Im}^{(2)},\R X)-b_3(\overline X_{\Im}^{{(2)}})$.
Assume that $\a_{\Im}^{(2)}\ge1$, then the universal coefficients
formula combined with (7-2) implies that $b_3(\overline
X_{\Im}^{(2)},\R X)\le b_3(\overline X_{\Im}^{(2)},\R X;\Z/2)\le
b_3(\C X_{\Im}^{(2)};\Z/2)$ and the duality gives $b_3(\overline
X_{\Im}^{(2)})= b_1(\overline X_{\Im}^{(2)};\d \overline
X_{\Im}^{(2)})= b_1(\overline X)+\a_{\Im}^{(2)}-1$.

Consider an  intermediate resolution, $\C X^{\res}_{\Im}\to\C X$,
which resolves the imaginary $\Z/2$-essential singularities, and
let $\C X^{\res}\to \C X^{\res}_{\Im}$ resolves all the other
singularities. Viewing $\C X^{(2)}_{\Im}$ as a subset of $\C
X^{\res}_{\Im}$, we obtain from the homology sequence of $(\C
X^{\res}_{\Im},\C X_{\Im}^{(2)})$ that $$ b_3(\C
X_{\Im}^{(2)};\Z/2)\le b_3(\C
X^{\res}_{\Im};\Z/2)+2\a_{\Im}^{(2)}-1 $$
%%%%%%%%%%%
and thus $$ \delta(\C X')\le b_3(\C
X^{\res}_{\Im};\Z/2)-b_1(\overline X)+\a_{\Im}^{(2)} $$

Denote by $\C E^{(2)}$ the exceptional divisor of the resolution
$\C X^{\res}\to\C X^{\res}_{\Im}$ and let $\theta=\dim\ker(H_2(\C
E^{(2)};\Z/2)\to H_2(\C X^{\res};\Z/2))$. The homology sequence of
$(\C X^{\res},\C E^{(2)})$ gives $$ b_3(\C
X^{\res}_{\Im};\Z/2)=b_3(\C X^{\res},\C E^{(2)};\Z/2)= b_3(\C
X^{\res};\Z/2)+\theta $$ Analyzing the homology sequence of $(\C
U_{\Re},\d\,\C U_{\Re})$, where $\C U_{\Re}$ is the union of the
components of $\C U$ around the real singularities of $\C X$, we
obtain an estimate $\theta\le d^{(2)}_{\Re}$, (here we use that
$\Z/2$-essential singularities of $\C X^{\res}_{\Im}$ are real).
Thus, $$ \delta(\C X')\le b_3(\C X^{\res};\Z/2)-b_1(\overline
X)+\a_{\Im}^{(2)}+ d^{(2)}_{\Re} \tag7-11 $$ In the case
$\a_{\Im}^{(2)}=0$, the same arguments bring an analogous formula,
in which $1$ stands instead of $\a_{\Im}^{(2)}$.
%$\delta(\C X')\le b_3(\C X^{\res}_{\Im};\Z/2)-b_1(\overline X)+1$).
This gives (4-24), because $b_1(\overline X)=\frac12b_1(\C X)$
(cf. 7.3). (4-25) follows from (7-11), (4-21) and the remark in
the beginning of 7.3.
%$$b_3(\overline X)=\frac12b_3(\C X)=\frac12b_3(\C X^{\res})=
%\frac12b_1(\C X^{\res})=h^{0,1}(\C X^{\res})$$
%%%%%%%%%%%%%%%%%%%%%%%%%%%

\subheading{7.5. Proof of Proposition 6.1.3} Along with the
arrangement, $\C A=\C A_1\cup\dots\cup\C A_m$, of hyperplanes $\C
A_i\subset\Cp{d}$ in generic position and the double covering
$X\to\Cp{d}$ branched along $\C A$, considered in 6.1, we consider
the affine arrangement $\C A_a=\C A\setminus\C A_m$ in
$\C^d=\Cp{d}-\C A_m$ and the double covering $X_a\to\C^d$ branched
along $\C A_a$. Endowed with one of the two real structures
covering the real structure in $\C^d$, $X_a$ becomes a real affine
variety, $\C X_a^\e\subset\C X^\e$, where $\e=+$ or $\e=1$. Put $$
\align \delta(\C X_a^\e)&=\dim\ker(H_d(\R X_a^\e;\Q)\to
H_d(X_a;\Q))\\ \delta_2(\C X_a^\e)&=\dim\ker(H_d(\R
X_a^\e;\Z/2)\to H_d(X_a;\Z/2))\\
\endalign
$$
%%%%%%%%%%

\tm{7.5.1. Lemma} $\delta(\C X_a)=\delta_2(\C X_a)=0$
\endtm

\demo{Proof} We obviously have
%The exact sequence of pair $(X_a,\R X_a^\e)$ implies that
$$ \align \delta(\C X_a^\e)&= \dim\ker(H_d(\R X_a^\e;\Q)\to
H_d(\overline X_a^\e;\Q))= b_{d+1}(\overline X_a^\e,\R
X_a^\e)-b_{d+1}(\overline X_a^\e)\\ \delta_2(\C X_a^\e)&\le
\dim\ker(H_d(\R X_a^\e;\Z/2)\to H_d(\overline X_a^\e;\Z/2))=
b_{d+1}(\overline X_a^\e,\R X_a^\e;\Z/2)-b_{d+1}(\overline
X_a^\e;\Z/2)\\
\endalign
$$ where $b_{d+1}(\overline X_a^\e,\R X_a^\e)\le b_{d+1}(\overline
X_a^\e,\R X_a^\e;\Z/2)$ by the universal coefficients formula. So,
it is enough to prove that $b_{d+1}(\overline X_a^\e,\R
X_a^\e;\Z/2)$ vanishes. Vanishing follows from the inequalities $$
\multline b_{d+1}(\overline X_a^\e, \R
X_a^\e;\Z/2)\le\sum_{k\ge1}b_{d+k}(X_a;\Z/2)
\le2\sum_{k\ge1}b_{d+k}(\C^d;\Z/2)+ \sum_{k\ge0}b_{d+k}(\C
A_a;\Z/2)=0
\endmultline
\tag7-13 $$ The first of these inequalities follows from (7-2)
applied to the complex conjugation, the second one follows from
(7-4) applied to the covering transform of $X_a\to\C^d$, and
$b_{d+k}(\C A_a;\Z/2)=0$ for $k\ge 0$ because $\R A_a$ is a
deformational retract of $\C A_a$. \qed\enddemo

%%%%%%%%%%%%%%%%%%%%%%%%%%%%%%%%

To derive Lemma 6.1.2 we note that the inclusion homomorphism
$H_d(X_a)\to H_d(X)$ is monomorphic, because $H_{d+1}(X,X_a)\cong
H^{d-1}(\C A_m)$ by the Alexander duality in $\Q$-homology
manifold $X$ and $H^{d-1}(\C A_m)=H^{d-1}(\Cp{d})=0$ for even $d$.
Vanishing of $\delta(\C X_a^\e)$ implies linear independence of
the fundamental classes, $[\G_i]\in H_d(\C X^\e)$, of
$\G_i=\pi^{-1}(W_i)$, for those of the partition components,
$W_i\subset\Rp{d}_\e$, which do not have points in common with $\R
A_m$. So, $\delta(\C X^\e)$ cannot exceed the number of the other
components, $W_i\subset\Rp{d}_\e$, such that $W_i\cap\R
A_m\ne\oo$. This number is obviously $N_{d-1}^{m-1}$, for a
generic hyperplane arrangement.

\rk{Remark} The most of the above arguments can be applied (after
a simple modification) to non-generic arrangement as well. The
exception is the argument in the last paragraph, which uses the
duality in $X$ and thus requires that $\C A$ is a
\sQQ-hypersurface.
\endrk

\heading \S8. Appendix B: Application of the Hodge theory to real
algebraic varieties
\endheading
\subheading{8.1. The Hodge structure in real \QR-varieties} For
\QR-varieties the mixed Hodge structure in $H^{*}(\C X;\C)$ is
known to be pure, which gives the usual Hodge splitting $H^*(\C
X;\C)=\oplus_{p,q\ge0}H^{p,q}(\C X)$ (indeed, the non-pure part of
$H^*(\C X;\C)$ is the kernel of $\res^*\:H^*(\C X\;\Q)\to H^*(\C
X^{\res};\Q)$ induced by a resolution $\res\:\C X^{\res}\to \C X$;
the latter is a degree $1$ map of rational homology manifolds, and
thus $\res^*$ is a monomorphism). An embedding $\C X\subset\Cp{N}$
gives a hyperplane class, $h\in H^{1,1}(\C X)$, and defines a
Lefschetz decomposition like in non-singular case. To define the
subspaces $P^{p,q}(\C X)\subset H^{p,q}(\C X)$ of primitive
$(p,q)$-classes in $\C X$, one can use the inner product in
$H^*(\C X;\C)$ induced from $H^*(\C X^{\res};\C)$ via $\res^*$.
The Hodge index theorem for $\C X^{\res}$ obviously descends to
$\C X$, since $\res^*$ preserves the Hodge filtration.

This allow us to can reproduce the arguments applied in \cite{Kh3}
to the case of non-singular real varieties as follows. Since
$\conj^*$ interchanges $H^{p,q}(\C X)$ and $H^{q,p}(\C X)$, both
the trace and the signature of the involution $\conj$, vanishes
being restricted to
 $H^{p,q}(\C X)\oplus H^{q,p}(\C X)$, for $p\ne q$. Therefore,
$$ \align \tr|_{H^*(\C X;\R)}&=\sum_{p=0}^{2n}\tr|_{H^{p,p}(\C
X)}\\ \s_{\conj}|_{H^{2n}(\C X;\R)}&=\s_{\conj}|_{H^{n,n}(\C X)}
\endalign
$$ Using that  multiplication by $h$, $\CD H^{p,p}(\C X)>\cup h>>
H^{p+1,p+1}(\C X)\endCD$, maps $\pm1$-eigenspace of $\conj^*$ into
the $\mp1$-eigenspace, we obtain $$ \align \tr|_{H^*(\C X;\R)}=
&\sum_{p=0}^{n}\tr|_{P^{p,p}(\C X)}\tag8-1\\ \tr|_{H^*(\C X;\R)}=
&\tr|_{H^{d}(\C X;\R)}+2\sum_{p\ge0}\tr|_{P^{n-1-2p,n-1-2p}(\C X)}
\tag8-2\\ \s_{\conj}|_{H^{2n}(\C X;\R)}=
&\sum_{p=0}^n\s_{\conj}|_{P^{p,p}(\C X)\wedge h^{n-p}}=
(-1)^n\sum_{p=0}^{n}\tr|_{P^{p,p}(\C X)} \tag8-3
\endalign
$$ Here $\s_{\conj}|_A$ denotes the signature of the involution,
$\s(A^+)-\s(A^-)$, for a $\conj^*$-invariant subspace $A\subset
H^d(\C X;\C)$ and the $\pm1$-eigenspaces, $A^\pm$, of $\conj^*$ in
$A$. To obtain (8-3) we used that the intersection form in
$H^{2n}(\C X;\R)$ is positive on $P^{p,q}(\C X)\wedge
h^{n-\frac12(p+q)}$ for even $\frac12(p+q)\le n$ and negative for
odd, which follows from the analogous fact (the Hodge index
theorem) for $\C X^{\res}$.

It is convenient to formulate the above identities in terms of the
functions $T^{\pm}(Z)=\chi(Z)\pm\s(Z)$, defined on compact
rational homology manifolds, $Z$, of dimension $4n$, (possibly,
with $\d Z\ne\oo$), and functions $D^\pm(\C X)=T^\pm(\overline
X)-\frac12T^\pm(\C X)$, defined on real \QR-varieties of dimension
$d=2n$, and on $\conj$-invariant compact codimension $0$
($\Q$-homology-)submanifolds of such varieties.

\tm{8.1.1. Theorem} Assume that $\C X$ is a real algebraic
\QR-variety of dimension $d=2n$, $\ee=(-1)^n$. Then $$ \align
D^{-\ee}(\C X)&=0\tag8-4\\ D^{\ee}(\C X)&=\chi(\R X)\tag8-5\\
b_d^{-\ee}(\overline X)&=\frac12(b_d^{-\ee}(\C X)-t(\C
X))\tag8-6\\ b_d^{\ee}(\overline X)&=\frac12(b_d^{\ee}(\C
X)+\chi(\R X)-t(\C X)) \tag8-7\\
\endalign
$$ where $t(\C X)=\sum_{p\ge0}\tr|_{P^{n-1-2p,n-1-2p}(\C X)}$.
\endtm

\demo{Proof of Theorem 8.1.1} (8-4) follows from (8-1) and (8-3).
Together with the Lefschetz fixed point formula for the involution
$\conj^*$, which can be stated as $D^{-\ee}(\C X)+D^{\ee}(\C
X)=\chi(\R X)$, it implies (8-5). (8-6) and (8-7) are versions of
(8-4) and (8-5), which are obtained by comparing (8-2) with (8-3).
\qed\enddemo

The Lefschetz hyperplane section theorem implies that $t(\C
X)=t(\Cp{d})=\kapnew$, if $\C X$ is a complete intersection of
dimension $d=2n$. This gives the relations (3-3).

\subheading{8.2. Proof of Theorem 3.2.2} One of the approaches to
prove Theorem 3.2.2 is to use a version of the arguments in 8.1
applied to the mixed Hodge structure in the Milnor fiber, $\C
U^\tau$ (using a similarity between the logarithm of the unipotent
part of the monodromy in $H^d(\C U^\tau;\C)$ and multiplication by
$h$ in $H^*(\C X;\C)$)). It is more elementary, however, to
present another proof, which uses intersection (of the middle
perversity) homology, rather then Mixed Hodge structure. Recall
that the intersection homology groups of algebraic varieties have
a pure Hodge structure, which is functorial, satisfy the usual
properties (the K\"ahler package), including the both Lefschetz
theorems, and Hodge index theorem \cite{Sa}, and coincides with
the usual Hodge structure if $\C X$ is a $\Q$-homology manifold.
This enables us to go through all the arguments of the previous
subsection and obtain a version of Theorem 8.1.1 for the
intersection homologies. For instance, a version of (8-4), that is
needed for the proof of Theorem 3.2.2, can be formulated as
follows. Denote by $\iota\chi$, $\iota\s$, the (middle perversion)
intersection homology  Euler characteristic and the signature of a
pseudo-manifold and put $ID^\pm(\C X)=IT^\pm(\C
X/\conj)-\frac12IT^\pm(\C X)$, where
$IT^{\pm}(Z)=\iota\chi(Z)\pm\iota\s(Z)$. Here $IT^{\pm}$ is
well-defined for any stratified pseudo-manifold, $Z$, of dimension
$4n$, whereas $\C X$ is supposed to be a real variety of dimension
$\dim_{\C}(\C X)=2n$, or $\conj$-invariant compact codimension $0$
(pseudo-)submanifold of such a variety.

\tm{8.2.1. Theorem} Assume that $\C X$ is a real algebraic variety
of dimension $d=2n$, $\ee=(-1)^n$. Then $ID^{-\ee}(\C X)=0$.
\qed\endtm

Assume now that $\C X$ is a real complete intersection which has
only isolated singularities. Let $\C X^\tau$ denote a real
deformation of $\C X$ and $\C U$, $\overline U$, $\C U^\tau$ and
$\overline U^\tau$ are chosen like in 3.3. By (8-4) and Theorem
8.2.1, $ID^{-\ee}(\C X)-D^{-\ee}(\C X^\tau)=0$. On the other hand,
additivity of $ID^{-\ee}$ and $D^{-\ee}$ together with
$D^{-\ee}(\C X-\C U)=ID^{-\ee}(\C X-\C U)$ imply that $D^{-\ee}(\C
U^\tau)-ID^{-\ee}(\C U)=0$, or equivalently, $$ T^{-\ee}(\C
U^\tau)-IT^{-\ee}(\C U)=2(T^{-\ee}(\overline U^\tau)
-IT^{-\ee}(\overline U)). $$ Furthermore, we have $$ \align
T^{-\ee}(\C U^\tau)&=b_0(\C U^\tau)+2b_d^{-\ee}(\C
U^\tau)+b_d^0(\C U^\tau) \tag8-8\\ IT^{-\ee}(\C U)&=b_0(\C
U^\tau)-b_{d-1}(\d\, \C U)\tag8-9\\
\endalign
$$ where (8-8) follows from that the connected component of $\C
U^\tau$ (homeomorphic to the corresponding Milnor fibers of $\C
X$) are homotopy equivalent to wedges of $d$-spheres. (8-9)
follows from vanishing of the intersection ``Betti numbers'',
$Ib_k(\C U)$, for $k\ge d$, from that $Ib_k(\C U)=b_k(\d\,\C U)$
for $k\le d$, \cite{GM1, p. 209}, and from vanishing of $\til
b_k(\d\, \C U)$ for all $k<d-1$ for ICIS.

The exact homology sequence of the pair $(\C U^\tau,\d\,\C
U^\tau)$ implies that $b_d^0(\C U^\tau)=b_{d-1}(\d\,\C U^\tau)$,
which gives $T^{-\ee}(\C U^\tau)-IT^{-\ee}(\C U)=2(b_d^{-\ee}(\C
U^\tau) +b_d^0(\C U^\tau))$,  since $\d\, \C U^\tau\cong \d\, \C
U$.
The same arguments can be applied to $\overline U$, due to Lemma
3.2.1; thus $T^{-\ee}(\overline U^\tau)-IT^{-\ee}(\overline U)
=2b_d^{-\ee}(\overline U^\tau)+2b_d^0(\overline U^\tau)$, which
proves (3-12) in Theorem 3.2.2. (3-13) follows from (3-12)
combined with the Lefschetz formula for the involution $\conj|_{\C
U^\tau}$.

\enddocument

\ref\key \by \paper \jour \vol \issue \yr \pages
\endref

\enddocument